\documentclass[11pt]{article}
\usepackage{amssymb}
\usepackage{mathrsfs}
\usepackage{latexsym,amsmath,amssymb,amsfonts,epsfig,graphicx,cite,psfrag}
\usepackage{eepic,color,colordvi,amscd}
\usepackage{ebezier}
\usepackage{verbatim}

\usepackage{subfigure}

\newcommand{\pf}{\noindent {\bf Proof.} }

\newtheorem{theorem}{Theorem}[section]
\newtheorem{lemma}[theorem]{Lemma}

\newtheorem{definition}[theorem]{Definition}

\definecolor{orange}{rgb}{1.0, 0.55, 0.0}

\def\qed{\hfill \rule{4pt}{7pt}}

\begin{document}

\title{Co-degree threshold for rainbow perfect matchings in uniform hypergraphs}

\author{Hongliang Lu\footnote{Partially supported by the National Natural
Science Foundation of China under grant No.11871391 and
Fundamental Research Funds for the Central Universities}\\
School of Mathematics and Statistics\\
Xi'an Jiaotong University\\
Xi'an, Shaanxi 710049, China\\
\medskip \\
Yan Wang\\
School of Mathematical Sciences\\
Shanghai Jiao Tong University\\
Shanghai 200240, China\\
\medskip \\
Xingxing Yu\footnote{Partially supported by NSF grant DMS-1954134}\\
School of Mathematics\\
Georgia Institute of Technology\\
Atlanta, GA 30332}

\date{}

\maketitle

\date{}

\maketitle

\begin{abstract}
Let $k$ and $n$ be two integers, with $k\geq 3$, $n\equiv 0\pmod k$, and $n$ sufficiently large.
 We
 determine the $(k-1)$-degree threshold for the existence of a rainbow
 perfect matchings in $n$-vertex
$k$-uniform hypergraph. This implies the result of R\"odl,
Ruci\'nski, and Szemer\'edi on the
$(k-1)$-degree threshold for the existence of  perfect matchings in
$n$-vertex $k$-uniform hypergraphs. In our proof, we identify the extremal configurations  of closeness, and consider whether or not the hypergraph is close to the extremal configuration. In addition, we also develop a novel absorbing device and generalize the absorbing lemma of R\"odl, Ruci\'nski, and Szemer\'edi.

\end{abstract}

\section{Introduction}

 A {\it hypergraph} $H$
consists of a vertex set $V(H)$ and an edge set $E(H)$ whose members
are subsets of $V(H)$. Let $H_1$ and $H_2$ be two hypergraphs. If
$V(H_1)\subseteq V(H_2)$ and $E(H_1)\subseteq E(H_2)$, then $H_1$ is
said to be a \emph{subgraph} of $H_2$ and we denote this by $H_1\subseteq H_2$.
Let $k$ be a positive integer and write $[k] = \{1,\ldots,k\}$.
For a set $S$, let ${S\choose k}:=\{T\subseteq S: |T|=k\}$.
A hypergraph $H$ is {\it $k$-uniform} if $E(H)\subseteq {V(H)\choose k}$, and a $k$-uniform hypergraph is also
called a {\it $k$-graph}. Given $T\subseteq V(H)$, let $H-T$ denote the subgraph of
$H$ with vertex set $V(H)\setminus T$ and edge set $\{e\in
E(H)\ : \ e\subseteq V(H)\setminus T\}$.

Let $H$ be a hypergraph and $S\subseteq V(H)$.
The {\it neighborhood} of $S$ in $H$ is $N_H(S)=\{T \subseteq
V(H)\setminus S:  T
\cup S \in E(H)\}$ and the {\it degree} of $S$ in $H$ is $d_H(S)=|N_H(S)|$.
For any postive integer $l$, $\delta_l(H):=\min\{d_H(S): S\in
{V(H)\choose l}\}$ is the {\it minimum $l$-degree} of $H$. Note that
$\delta_1(H)$ is called the minimum \emph{vertex degree} of $H$. If $H$ is a
$k$-graph then $\delta_{k-1}(H)$ is known as the minimum
\emph{co-degree}  of $H$. For a subset $M\subseteq E(H)$, we let
$V(M)=\bigcup_{e\in M}e$.

A \emph{matching} in a hypergraph $H$ is a subset of $E(H)$ consisting
of pairwise disjoint edges, which  is   \emph{perfect}
if $V(M)=V(H)$. While a maximum matching in a graph can be found in
polynomial time \cite{E65}, it is NP-hard to find even for 3-graphs
\cite{Ka72}. Much effort has been devoted to finding good sufficient
conditions for the existence of a large matching in uniform
hypergraphs, including Dirac type conditions. A celebrated result in
this area is due to
R\" odl, Ruci\'nski, and Szemer\'edi \cite{RRS09}, which refines the analysis in \cite{RRS06}. They determined the
minimum co-degree threshold function that ensures a perfect matching
in $n$-vertex $k$-graphs. For integers $k,n$, with $n\geq k\ge 3$ and $n\equiv
0\pmod k$, let
\[
t(n,k):=\begin{cases}
n/2 + 2 - k, & \text{if $k/2$ is even and $n/k$ is odd,}\\
n/2 + 3/2 - k, & \text{if $k$ is odd and $(n-1)/2$ is odd,}\\
n/2 + 1/2 - k, & \text{if $k$ is odd and $(n-1)/2$ is even,}\\
n/2 + 1 - k, & \text{otherwise.}
\end{cases}
\]
R\" odl, Ruci\'nski, and Szemer\'edi \cite{RRS09} proved the following
result.

\begin{theorem} [R\" odl, Ruci\'nski, and Szemer\'edi 2009] \label{RRS}
Let $k,n$ be integers, with $k\ge 3$, $n\equiv
0\pmod k$, and $n$ sufficiently large. Let $H$ be a $k$-graph on $n$
vertices such that $\delta_{k-1}(H)> t(n,k)$. Then $H$ has a perfect matching.
\end{theorem}

Codegree condition $\delta_{k-1}(H)> t(n,k)$ is best possible because of the following
$k$-graphs $H(n,k)$ on vertex set $[n]$ from \cite{Kuh06} (for odd $k$) and \cite{RRS09}
(for even $k$): When $k$ is odd, $[n]$ has a partition $A,B$
such that $|A|$ is the unique odd integer from the set
$\{\frac{n-2}{2},\frac{n-1}{2},\frac{n}{2},\frac{n+1}{2}\}$
and $E(H(n,k))=\left\{e\in {V(H(n,k))\choose k} : |e\cap A| \equiv 0 \pmod
2\right\}$.  When $k$ is even, $V(H(n,k))$ has a partition $A,B$ such that
$$|A|:=\begin{cases}
n/2-1,& \text{if  $n/k$ is even,}\\
n/2-1,& \text{if $n/k$ is odd and $n/2$ is odd,}\\
n/2,& \text{if $n/k$ is odd and $n/2$ is even,}
\end{cases}$$
and  $E(H(n,k))=\{e\in {V(H(n,k))\choose k} : |e\cap A| \equiv 1 \pmod
2\}$. Note that the sets $A,B$ are called {\it partition classes} of
$H(n,k)$.

Let $\mathcal{F} = \{F_1,\ldots, F_t\}$ be a family
of hypergraphs; a set of pairwise disjoint edges, one from
each $F_i$, is called a \emph{rainbow matching} for $\mathcal{F}$. In
this case, we also say that ${\cal F}$ {\it admits} a rainbow
matching. There has been effort to extend results on matchings in hypergraphs to rainbow
matchings, see for instance, \cite{AH17, FK20, HLS12, GLMY21, KK21, KLLM19, KLLM21, KL17, Kup21, Lu20, Lu21, LYY}.
The main result in this paper is a rainbow version of
Theorem~\ref{RRS}.


\begin{theorem}\label{main}
Let $k,n$ be integers with $k\ge 3$, $n\equiv 0\pmod k$, and $n$
sufficiently large.  Let $\{F_1,\ldots, F_{n/k}\}$ be a family of
$k$-graphs on the common vertex set $[n]$, such that
$\delta_{k-1}(F_i) >t(n,k)$ for $i\in [n/k]$.
Then $\{F_1,\ldots, F_{n/k}\}$ admits a rainbow perfect matching.

\end{theorem}

It is easy to see that we derive Theorem~\ref{RRS} from
Theorem~\ref{main} by setting $F_1=\ldots =F_{n/k}=H$.  Moreover, if $F_i=H(n,k)$ for $i\in [n/k]$ then $\{F_1, \ldots,
F_{n/k}\}$ admits no rainbow perfect matching.  So the co-degree bound
in Theorem~\ref{main} is best possible. We point out that
Theorem~\ref{main} for $k=2$ is a result of Joos and  Kim
\cite{JK20} and Akiyama and Frankl \cite{AF85}.

For $n\equiv 0\pmod k$, let $\{F_1,\ldots,
F_{n/k}\}$ be a family of $k$-graphs on the same vertex set $[n]$.
Let $X=\{x_1, \ldots, x_{n/k}\}$ be disjoint from $[n]$.
We consider the hypergraph ${\cal F}(n,k)$
with vertex set $X\cup [n]$ and edge set $\bigcup_{i=1}^{n/k}
\{\{x_i\}\cup e\ :\ e\in E(F_i)\}$. We denote this hypergraph by
${\cal H}(n,k)$ when $F_i=H(n,k)$ for $i\in [n/k]$ with same partition classes
$A,B$ of $[n]$, and refer to ${\cal H}(n,k)$ as {\it extremal
  configuration}.  It is easy to see the following is true.

\medskip

\textbf{Observation}. $\delta_{k-1}(F_i)>t(n,k)$ for $i\in [n/k]$ implies that $d_{\cal F}(n,k)(S)>t(n,k)$ for any $S\in {V({\cal F}(n,k))\choose k}$ with $|S\cap X|=1$.  $\{F_1,\ldots,F_{n/k}\}$ admits a rainbow matching if, and only if,  $\mathcal{F}(n,k)$ has a perfect matching.

So we will show that ${\cal F}(n,k)$ has a perfect
matching. Indeed, we consider a more general class of hypergraphs.  Let $Q,V$ be two disjoint sets.  A $(k+1)$-graph $H$ with vertex set $Q \cup V$
is said to be \emph{$(1,k)$-partite} with {\it partition classes} $Q,V$ if, for each edge $e\in E(H)$, $|e\cap Q|=1$ and $|e\cap V|=k$.
A $(1,k)$-partite $(k+1)$-graph $H$ with partition classes $Q,V$ is  \emph{balanced} if $|V|=k|Q|$.
We say that a subset $S \subseteq V(H)$ is \textit{balanced} if
$|S\cap V|=k|S\cap Q|$. Theorem~\ref{main} follows from the following
result.

\begin{theorem}\label{general}
Let $k,n$ be integers with $k\ge 3$, $n\equiv 0\pmod k$, and $n$
sufficiently large. Let ${\cal F}$ be a balanced $(1,k)$-partite $(k+1)$-graph with
partition classes $X, [n]$, such that for any $S\in {V({\cal F})\choose
  k}$ with $|S\cap X|=1$, $d_{{\cal F}}(S)>t(n,k)$. Then ${\cal
  F}$ admits a perfect matching.
\end{theorem}

 To prove Theorem~\ref{general}, we consider whether or not ${\cal F}$
 is ``close'' to the extremal configuration
${\cal H}(n,k)$. In Section 2, we describe several properties of the
extremal configurations. In Section 3, we prove Theorem~\ref{general}
for the case when ${\cal F}$ is close to the extremal
configuration. In Section 4, we study  absorbing devices for perfect
matchings, and in Section 5, we study an absorbing device for near
perfect matchings.
We deal with the case when ${\cal F}$ is not close to the extremal configuration in Section 6 and offer some concluding remarks  in Section 7.

\section{Properties of Extremal configurations}

We will often use the following $(1,k)$-partitie $(k+1)$-graphs as
intermediate configuration to compare  $(1,k)$-partitie $(k+1)$-graphs with ${\cal H}(n,k)$. Suppose $W,U$ form a partition of $[n]$ such that $|W|= (1/2 \pm
o(1))n$ and $|U|=(1/2\pm o(1)) n$. For $i\in \{0,1\}$, let
 $H_{n,k}^{i}(W,U)$ denote the $k$-graph with vertex set $[n]$ and
 edge set $\{S\in {[n]\choose k}: |S\cap W|\equiv i \pmod 2\}$.
When $|W|=\lfloor n/2\rfloor$, 
we denote
$H_{n,k}^{i}(W,U)$  by $H_{n,k}^{i}$.

We need the following definition to quantify the difference between ${\cal
  F}(n,k)$ and ${\cal H}(n,k)$.
Let $\varepsilon > 0$ be  a real number. Given two $k$-graphs $H_1, H_2$ with $V(H_1)=V(H_2)$, we say that $H_2$ is \textit{strongly $\varepsilon$-close} to $H_1$  if $|E(H_1)\setminus E(H_2)|\leq \varepsilon |V(H_1)|^k$.
We say that $H_2$ is \textit{weakly $\varepsilon$-close} to $H_1$ if
$c(H_1,H_2)\leq \varepsilon |V(H_1)|^k$, where $c(H_1,H_2)$ be the minimum of $|E(H_1)\backslash E(H_2')|$
taken over all isomorphic copies $H_2'$ of $H_2$ with  $V(H_2') =
V(H_2)$. It is easy to see that the following is true.


\begin{lemma}
Let $\varepsilon > 0$ be a real number.
Let $k,n$ be integers with $k\ge 3$, $n\equiv 0\pmod k$ and $n$ is sufficiently large.
Let
$W,U$ be a partition of $[n]$ with $|W|= (1/2\pm o(1))n$ and $|U|=(1/2\pm o(1))n$.
Then the following
statements hold.
\begin{itemize}
\item [(i)] If  $k$ is odd then $H_{n,k}^{i}(U,W)$ and $H_{n,k}^{j}(W,
  U)$, $i,j\in \{0,1\}$,  are weakly $\varepsilon$-close to each
  other.
\item [(ii)] If $k$ is even then $H_{n,k}^{i}(U,W)=H_{n,k}^{i}(W,U)$
  for $i\in \{0,1\}$.
\end{itemize}
\end{lemma}

Let $k, n$ be integers with $k\ge 3$ and $n\equiv 0 \pmod k$.
Let $m_0,m_1$ be integers between $0$ and $n/k$ (inclusive).
For convenience,  define $\mathcal{H}_{n,k}^{m_0,m_1}(W,U)$ as the
$(1,k)$-partite $(k+1)$-graph $\mathcal{H}$ with partition classes $X$ and $[n]$ and
a partition $W,U$ of $[n]$, such that $|X| = n/k$ and $|\{x\in X:
N_{\mathcal{H}}(x) = H_{n,k}^i(W,U)\}|=m_i$ for $i\in \{0,1\}$. For
$i\in \{0,1\}$, if $m_i=m$ and $m_0+m_1=n/k$ then we denote
$\mathcal{H}_{n,k}^{m_0,m_1}(W,U)$ by $\mathcal{H}_{n,k}^{i}(W,U;
m)$.

In the remainder of this section, we study 
$(1,k)$-partite $(k+1)$-graphs that ${\cal F}$ are close to some ${\cal H}(n,k)$.
and consider those vertices in ${\cal F}$ that is
contained in lots of edges of ${\cal H}(n,k)$.   So we introduce the
following concept.
Let $k,n$ be integers with $k \ge 3$, $n \equiv 0 \pmod k$, and $n$
sufficiently large.
Let $\mathcal{F}$ and ${\cal H}$ be $(1,k)$-partite $(k+1)$-graphs
with partition classes $X, [n]$.
A vertex $v$ of ${\cal F}$ is said to be {\it $\alpha$-good with respect
to} ${\cal H}$ if
$| N_{{\cal H}}(v) \setminus  N_{{\cal F}}(v) | < \alpha |V(\mathcal{F})|^{k}$. 
Otherwise, $v$ is said to be {\it
  $\alpha$-bad} with respect to ${\cal H}$.


The following lemma shows that the number of bad vertices in
${\cal F}$ is small if ${\cal F}$ is close to
$\mathcal{H}^i_{n,k}(W,U;m)$ for some $i\in \{0,1\}$.

\medskip

\begin{lemma} \label{bad-vertices}
Let $k, n$ be integers with $k\ge 3$  and $n \equiv 0 \pmod k$, and let $\varepsilon$ be a constant such that $0<1/n\ll \varepsilon\ll 1/k$.
Let $\mathcal{F}$ be a $(1,k)$-partite $(k+1)$-graph with partition classes 
$X$ and $[n]$ where $|X|=n/k$. 
Let $0 \le m \le n/k$ be an integer.
If $\mathcal{F}$ is strongly $\varepsilon$-close to some
$\mathcal{H}^i_{n,k}(W,U;m)$, where $|W| = n/2 \pm o(n)$ and $|U|= n/2 \pm o(n)$,
then  the number of
$\varepsilon^{2/3}$-bad vertices in $\mathcal{F}$ with respect to
$\mathcal{H}^i_{n,k}(W,U;m)$   is at most $(1+1/k)(k+1)\varepsilon^{1/3}n$.
\end{lemma}

\pf
Let $N$ be the set of $\varepsilon^{2/3}$-bad vertices in ${\cal F}$ with respect to $\mathcal{H}^i_{n,k}(W,U;m)$.
If $|N|>((1+1/k)(k+1)\varepsilon^{1/3}n$ then
\begin{align*}
|E(\mathcal{H}^i_{n,k}(W,U;m)) \setminus E(\mathcal{F})| > \frac{1}{k+1}|N|\varepsilon^{2/3} |V(\mathcal{F})|^k \ge \varepsilon \left(n+\frac{n}{k}\right)^{k+1}.
\end{align*}
This  contradicts the assumption that  $\mathcal{F}$ is strongly $\varepsilon$-close to 
$\mathcal{H}^i_{n,k}(W,U;m)$.
\qed

\medskip

The next lemma says that we can find an edge in $\mathcal{F}$ which serves as ``parity breaker".
This is the only place in the proof of Theorem~\ref{general} where we require $\delta_{k-1}(N_{\mathcal{F}}(x))>t(n,k)$ for all $x \in X$.


\begin{lemma} \label{parity}
Let $k,n$ be integers with $k \ge 3$, $n \equiv 0
\pmod k$, and $n \ge 2k$. Let $\mathcal{F}$ be a balanced $(1,k)$-partite $(k+1)$-graph with partition
classes $X, [n]$, such that  $\delta_{k-1}(N_{\mathcal{F}}(x))>t(n,k)$ for all $x\in X$.
Then for any proper subset $X'$ of $X$  and any partition $W,U$ of
$[n]$ with $\min\{|W|,|U|\}\geq k$,
there exists $e_0$ with $e_0 \in E(\mathcal{F})$ or $e_0=\emptyset$,  such that, for $i\in \{0,1\}$,
\begin{itemize}
\item[(i)] if $i=0$ then $|W \setminus  e_0 | \equiv |X' \setminus e_0| \pmod 2$;
\item[(ii)] if $i=1$ then $|W \setminus  e_0| \equiv |X \setminus  (X' \cup e_0)| \pmod 2$;
\item[(iii)] if $k-i$ is even then $|U \setminus  e_0| \equiv |X' \setminus e_0| \pmod 2$;
\item[(iv)] if $k-i$ is odd then $|U \setminus  e_0| \equiv |X  \setminus (X' \cup e_0)| \pmod 2$.
\end{itemize}
\end{lemma}

\pf
Let $q := |X'|$ and fix  $x\in X\setminus X'$. Then
$\delta_{k-1}(N_{\mathcal{F}}(x)) \ge t(n,k) + 1$.
 In all cases below, we may assume the assertion of this lemma does
not hold for $e_0=\emptyset$.

\medskip

\textbf{Case 1}: $i=0$ and $k-i=k$ is even.

Then (i) and (iii) are relevant and not both true for $e_0 =
\emptyset$. Hence,  $|W| \not \equiv q \pmod 2$ or $|U| \not \equiv q \pmod 2$.
Therefore, since $n \equiv 0 \pmod 2$ (as $n \equiv 0 \pmod k$ and $k$
is even),
$|W| \not \equiv q \pmod 2$, $|U| \not \equiv q \pmod 2$ and $|W|\equiv |U| \pmod 2$.

If there exists $e_0 \in E({\cal F})$ with $x\in e_0$ such that
$|e_0 \cap W| =1$ or $|e_0 \cap U| = 1$, then both (i) and (iii)
holds with this $e_0$. So we may assume that, for every $e\in E({\cal
  F})$ with $x\in e$, $|e\cap W| \ne 1$ and $|e \cap U| \ne 1$.

Then  for any $(k-1)$-set $S\subseteq U$, $N_{N_{{\cal F}}(x)}(S)\subseteq
U\setminus S$. Thus,  $|U-S| \ge |N_{N_{\cal F}(x)}(S)| \ge t(n,k)
+1\ge n/2+2-k$. Thus $|U| \ge n/2 + 1$.

Similarly, we derive $|W| \ge n/2 + 1$. This leads to a contradiction as $n = |U| + |W|$.

\medskip

\textbf{Case 2}: $i=0$ and $k-i=k$ is odd.

Then  (i) and (iv) are relevant and not both true for $e_0 =
\emptyset$. So  $|W| \not \equiv q \pmod 2$ or $|U| \not \equiv n/k -
q \pmod 2$. Note that $n/k-q\equiv n-q\pmod 2$.
Thus since $|U| + |W| = n$, we have $|U| \not \equiv n - q \pmod 2$ and $|W| \not \equiv q \pmod 2$.

\medskip

\textit{Subcase 2.1}: $n$ is even; so  $n/k$ is even.

If there exists $e_0\in {\cal F}$ such that $x\in e_0$ and $|e_0 \cap
W| = 1$ or $|e_0 \cap W| = k - 2$ then (i) and (iv) holds for this
$e_0$. So assume that for any $e\in E({\cal F})$ with $x\in e$, $|e\cap
W| \ne 1$ and $|e\cap W| \ne  k - 2$.

Then  for any $(k-1)$-subset $S\subseteq U$, $N_{N_{\mathcal{F} }(x)}(S)\subseteq
U\setminus S$. Thus, $|U\setminus S| \ge
|N_{N_{{\cal F}}(x)}(S)| \ge t(n,k)+1 \ge n/2+2-k$; so $|U| \ge n/2 + 1$.

Similarly,  for any $(k-1)$-subset $T\subseteq [n]$ with $|T \cap W|
= k-2$, $N_{N_{{\cal F}}(x)}(T)\subseteq W\setminus T$. Thus, $|W\setminus T| \ge |N_{N_{{\cal F}}(x)}(T)| \ge t(n,k) +1\ge
n/2+2-k$; so $|W| \ge n/2$. However, this is a contradiction as $n = |W| + |U| \ge n + 1$.

\medskip

\textit{Subcase 2.2}. $n$ is odd; so $n/k$ is odd.

Suppose $q\equiv 0\pmod 2$. If there exists $e_0
\in E({\cal F})$  such that $x\in e_0$ and $|e_0 \cap W| = 1$ or $|e_0 \cap W| = k
- 2$ then (i) and (iv) hold for this $e_0$. So assume that for any
$e\in E({\cal F})$ with $x\in e$,  $|e \cap W| \ne 1$ and $|e \cap W| \neq k
- 2$.
Then, for any $(k-1)$-subset $S\subseteq U$, $N_{N_{{\cal
    F}}(x)}(S)\subseteq U\setminus S$ which implies $|U\setminus S| \ge |N_{N_{{\cal
      F}}(x)}(S)| \ge t(n,k) +1\ge n/2+3/2-k$; so $|U| \ge (n+1)/2$.
Similarly, for any  $(k-1)$-subset $T\subseteq [n]$ with $|T \cap W| =
k-2$, $N_{N_{{\cal F}}(x)}(T)\subseteq W\setminus T$ and, hence, $|W\setminus T| \ge |N_{N_{{\cal F}}(x)}(T)| \ge t(n,k)+1 \ge
n/2+3/2-k$; so $|W| \ge (n-1)/2$ and the equality holds only when
$(n-1)/2$ is even.
Thus since $n = |W| + |U|$,  $|U| = (n+1)/2$ and $|W| =
(n-1)/2\equiv 0 \pmod 2$, which leads to a contradiction since  $|W|\not \equiv q \pmod 2$ and  $q\equiv 0\pmod 2$.

Now assume  $q \equiv 1 \pmod 2$. Then we have $X' \ne \emptyset$, and
fix $x'\in X'$.
If there exists an edge $e_0\in E({\cal F})$ such that  $x'\in e_0$ and $|e_0
\cap W| = 2$ or $|e_0 \cap W| = k - 1$, then (i) and (iv) hold with
this $e_0$. So we may assume that for any $e\in E({\cal F})$ with
$x'\in e$,  $|e\cap W|\ne 2$ and $|e\cap W| \ne k - 1$.
Then for any $(k-1)$-subset $T\subseteq W$, $N_{N_{{\cal
F}}(x')}(T)\subseteq W\setminus T$ and, hence, $|W\setminus T| \ge
|N_{N_{{\cal F}}(x')}(T)| \ge t(n,k) +1\ge n/2+3/2-k$; so $|W| \ge (n+1)/2$.
Similarly, for any $(k-1)$-subset $S\subseteq U$ with $|S \cap U| =
k-2$, $N_{N_{{\cal F}}(x')}(S)\subseteq U\setminus S$ and, hence,  $|U\setminus S| \ge |N_{N_{\cal F}(x')}(S)| \ge t(n,k) +1\ge
n/2+3/2-k$; so $|U| \ge (n-1)/2$ and the equality holds only when $(n-1)/2$ is even.
Thus since $n = |W| + |U|$,  $|U| = (n-1)/2\equiv 0\pmod 2$ and $|W| = (n+1)/2$,
 which leads to a contradiction  since $|U|\not\equiv n-q\pmod 2$ and $n-q\equiv 0\pmod 2$.

\medskip

\textbf{Case 3}: $i=1$ and $k-i=k-1$ is odd.

So  (ii) and (iv) are relevant and not both true for $e_0 =
\emptyset$. Thus, $|W| \not \equiv n/k - q \pmod 2$ or $|U| \not
\equiv n/k - q \pmod 2$. Note that $n$ is even as $k$ is even.
Since $n = |W| + |U|$, $|W| \not \equiv n/k - q \pmod 2$, $|U| \not \equiv n/k - q \pmod 2$ and $|W|\equiv |U|\pmod 2$.

Suppose $q \equiv 1 \pmod 2$. Then $X' \ne \emptyset$, and fix $x'\in
X'$. If there exists  $e_0 \in E({\cal F})$ such that $x'\in e_0$, and
$|e_0 \cap W| = 1$ or $|e_0 \cap W| = k-1$, then (ii) and (iv) holds
for this $e_0$, in this case, $X\backslash(X'\cup e_0)=X\backslash X'$. So assume that,  for any
$e\in E({\cal F})$ with $x'\in e$,  $|e\cap W| \ne 1$ and $|e\cap
W| \ne  k-1$. Hence,  for any $(k-1)$-subset $S\subseteq [n]$ with
$S\subseteq W$, $N_{N_{{\cal F}}(x')}(S)\subseteq W\setminus S$ and,
hence, $|W\setminus S| \ge |N_{N_{{\cal F}}(x')}(S)| \ge
t(n,k) +1\ge n/2+2-k$; so $|W| \ge n/2 + 1$.
Similarly, $|U| \ge n/2 + 1$.
This is a contradiction since $n = |W| + |U| \ge n + 2$.

So $q\equiv 0\pmod 2$. If there exists $e_0 \in E({\cal F})$
with $x\in e_0$ such that $|e_0 \cap W| = 2$ or $|e_0 \cap U| = 2$
then (ii) and (iv) holds with this $e_0$. Hence, we may assume that, for any
$e\in E({\cal F})$ with $x\in e$,  $|e \cap W| \ne 2$ and $|e\cap
U| \ne  2$.

Hence,  for any $(k-1)$-subset $S\subseteq [n]$ with $|W\cap S|=k-2$,
$N_{N_{{\cal F}}(x)}(S)\subseteq W\setminus S$ and, hence,
$|W\setminus S| \ge |N_{N_{{\cal F}}(x)}(S)| \ge t(n,k) +1\ge
n/2+2-k$. So $|W| \ge n/2$ and  equality holds only when $n/k$ is
even or when $n/k$ is odd and $k/2$ is odd.

Similarly, $|U| \ge n/2$ and equality holds only when $n/k$ is even
or when $n/k$ is odd and $k/2$ is odd.
Since $n = |W| + |U|$, $|W| = |U| = n/2$.
If $n/k$ is even then $|W|=n/2$ is even (as $k$ is even); however,
$|W|$ is odd since $|W| \not \equiv n/k - q \pmod 2$, a
contradiction.  If $n/k$ is odd then $k/2$ is odd and, hence, $|W|=n/2$ is
odd. However, $|W|$ is even, since $|W| \not \equiv n/k - q \pmod 2$,  a contradiction.

\medskip

\textbf{Case 4}: $i=1$ and $k-i=k-1$ is even.

Then  (ii) and (iii) are relevant, and not both true for $e_0 =
\emptyset$. Hence, $|W| \not \equiv n/k - q \pmod 2$ or $|U| \not \equiv q \pmod 2$.
Note that $n = |W| + |U|$ and $n/k-q\equiv n-q \pmod 2$. So  $|W| \not
\equiv  n - q \pmod 2$ and $|U| \not \equiv q \pmod 2$.

\medskip

\textit{Subcase 4.1}: $n$ is even.

If there exists $e_0\in E({\cal F})$ with $x\in e_0$ such that  $|e_0
\cap U| = 1$ or $|e_0 \cap U| = k - 2$, then (ii) and (iii) holds
for this $e_0$. So assume that for any $e\in E({\cal F})$ with $x\in
e$,  $|e \cap U| \neq 1$ and $|e \cap U| \neq k - 2$.

Then, for any $(k-1)$-subset $S\subseteq W$, $N_{N_{{\cal
F}}(x)}(S)\subseteq W\setminus S$ and, hence,  $|W\setminus S| \ge
|N_{N_{{\cal F}}(x)}(S)| \ge t(n,k)+1 \ge n/2+2-k$; so $|W| \ge n/2 + 1$.
Similarly, for any $(k-1)$-subset $T\subseteq [n]$ with $|T \cap U| =
k-2$, $N_{N_{{\cal F}}(x)}(T)\subseteq U\setminus T$ and, hence,  $|U\setminus T| \ge |N_{N_{{\cal F}}(x)}(T)| \ge t(n,k) +1\ge
n/2+2-k$; so $|U| \ge n/2$. Now $n = |W| + |U| \ge n + 1$, a contradiction.

\medskip

\textit{Subcase 4.2}: $n$ is odd.

Suppose $q\equiv 0\pmod 2$. If there exists $e_0\in E({\cal F})$ with
$x\in e_0$ such that  $|e_0 \cap U| = 1$ or $|e_0 \cap U| = k-2$
then (ii) and (iii) hold for this $e_0$. So assume that for any $e\in
E({\cal F})$ with $x\in e$,  $|e \cap U| \ne 1$ and  $|e\cap U| \ne  k -2$.
Then for any some $(k-1)$-subset $S\subseteq W$, $N_{N_{\cal F}(x)}(S)\subseteq W\setminus S$ and, hence, $|W\setminus S| \ge
|N_{N_{{\cal F}}(x)}(S)| \ge t(n,k) +1\ge n/2+3/2-k$;
so $|W| \ge n/2 + 1/2$. Also, for any $(k-1)$-subset $T\subseteq [n]$
with $|T \cap U| = k-2$, $N_{N_{{\cal F}}(x)}(T)\subseteq U\setminus T$
and, hence,
 $|U\setminus T| \ge |N_{N_{{\cal F}}(x)}(T)| \ge t(n,k)+1 \ge n/2+3/2-k$; so
$|U| \ge n/2-1/2$ and the equality holds only when $(n-1)/2$ is even.
Since $n = |W| + |U|$,
$|W| = (n+1)/2$ and $|U| = (n-1)/2$ and $(n-1)/2$ is even.
However, this is a contradiction as $|U| \not \equiv q \pmod 2$.

So $q \equiv1\pmod 2$. Then $X' \ne \emptyset$, and fix $x'\in
X'$.  If there exists $e_0\in E({\cal F})$ with
$x'\in e_0$ such that $|e_0 \cap W| = 1$ or $|e_0 \cap W| = k - 2$,
then (ii) and (iii) hold for this $e_0$. So assume that for any $e\in
E({\cal F})$ with $x'\in e$, $|e \cap W| \ne  1$ and $|e_0 \cap W|
\ne k - 2$.

Then for any $(k-1)$-subset $S\subseteq U$, $N_{N_{{\cal F}}(x')}(S)\subseteq U\setminus S$ and, hence,  $|U\setminus S| \ge
|N_{N_{{\cal F}}(x')}(S)| \ge t(n,k) +1\ge n/2+3/2-k$; so $|U| \ge (n+1)/2$.
Similarly, for any $(k-1)$-subset $T\subseteq [n]$ with $|T \cap W| =
k-2$, $N_{N_{{\cal F}}(x')}(T)\subseteq W\setminus T$ and, hence,  $|W\setminus T| \ge |N_{N_{{\cal F}}(x')}(T)| \ge t(n,k)+1 \ge
n/2+3/2-k$; so $|W| \ge (n-1)/2$ and the equality holds only when $(n-1)/2$ is even.
Since $n = |W| + |U|$, we have $|W| = (n-1)/2\equiv 0\pmod 2$ and $|U| = (n+1)/2$, which is a contradiction since $|W|\neq n-q \pmod 2$.
\qed

\section{Hypergraphs close to extremal configurations}

We often need to move some vertices between two sets and keep track of their degrees.
The following notation will be convenient.
Let $H$ be a $k$-graph.
For $j \in \{0,1\}, v \in V(H)$, and $S \subseteq V(H)$, we define
$$d_{H,S}^j(v) := |\{ e \in E(H): v \in e \text{ and } |e \cap S| \equiv j \pmod 2\}|.$$

We begin with a lemma that allows us to find a matching in the
hypergraphs in question covering any
small fixed set of vertices. For convenience, we set the following
parameters for a  given integer $k$ for the remainder of this section:
$\eta = 1/(4k!)$,
$c = 1/(8(k+1)!)$,
$\varepsilon = 1/(80^kk^{k-5}(k!)^k((k+1)!)^k)^{3/2}$,
and $\gamma = \frac{\varepsilon^{2/3} 2^k}{c^k k^k}$.

\begin{lemma} \label{match-bad-vertices}
Let $k,n$ be integers with $k \ge 3$,  $n \equiv 0 \pmod k$, and $n \gg
1/\varepsilon$.
Let $\mathcal{F}$ be a balanced $(1,k)$-partite $(k+1)$-graph with partition
classes $X, [n]$,  and assume  $\delta_{k-1}(N_{\mathcal{F}}(x)) > t(n,k)$ for all $x \in X$.
Let $W, U$ be a partition of $[n]$ such that $\min\{|W|,|U|\}\geq k$,
and  let $X_0, X_1$ be a partition of $X$ such that $d_{\mathcal{F}, W}^j(x)\geq \eta n^{k}$ for $j\in \{0,1\}$ and $x\in X_j$.
Suppose there exists  $i\in \{0,1\}$ such that $d_{\mathcal{F} -
  X_{1-i}, W}^i(v)\geq \eta n^{k}$ for all $v\in [n]$.

Then for any $e_0$ satisfying the conclusion  of Lemma~\ref{parity}
and for any $N \subseteq X \cup [n]$ with $|N|\le 2c n$,
there exists a matching $M$ in $\mathcal{F}$ such that $N
\subseteq V(M)$, $e_0\in M$ when $e_0\ne \emptyset$,
$|V(M)| \le (k+1) 2c n$, and $|W \setminus  V(M) |\equiv  |X_{1} \setminus V(M)| \pmod 2.$
\end{lemma}

\pf
Let $N_{0}=N \cap X_0$ and $N_{1}=N \cap X_1$. 
Let $e_0$ be the edge satisfying the conclusion of Lemma~\ref{parity}.
If $N \setminus e_0 = \emptyset$ then let $M := \emptyset$ when
$e_0=\emptyset$, and  $M=\{e_0\}$ when $e_0\ne \emptyset$. So assume
$N\setminus e_0\ne \emptyset$.
Divide $N\setminus
e_0$ into three pairwise disjoint sets: $N_{0} \setminus e_0=\{v_1,\ldots,v_r\}$,
$N_{1}\setminus e_0=\{v_{r+1},\ldots,v_{s}\}$, and $(N\setminus
e_0)\cap [n] = \{v_{s+1},\ldots,v_{t}\}$. We find the desired matching  $M$ by
covering the vertices in $N\setminus e_0$ greedily.

Suppose we have found the matching $M_u:=\{e_0,e_1,...,e_u\}$ for some
$u \ge 0$ such that for $1\leq j\leq u$,  $\{v_j\}= e_j\cap N$ and  for $1\leq i\leq \min\{u, r\}$,   $|e_j\cap W|\equiv 0 \pmod 2$, for $r+1 \le j \leq \min\{u,s\}$, $|e_j\cap W|\equiv 1 \pmod 2$, and for $s+1 \le j \leq \min\{u,t\}$, $|e_j\cap W|\equiv i \pmod 2$ and $e_j\cap X\subseteq X_i$.

If $N \subseteq \cup_{i=0}^u e_i$, then $M_u$ is the desired matching.
Otherwise, let $v_{u+1} \in N \setminus \cup_{i=0}^u e_i$.
Since $u+1<|N| \le 2c n$, the number of edges in $\mathcal{F}$ containing $v_{u+1}\in N$ and a vertex from $\cup_{i=0}^u e_i\cup N$ is less than
$$\left(|N|+(u+1)(k+1)\right)n^{k-1} \le (k+1) 2c n^k.$$

By assumption, $d_{\mathcal{F}, W}^0(v)\geq \eta n^{k}$ for $v\in N_{0}$,
$d_{\mathcal{F}, W}^1(v)\geq \eta n^{k}$ for $v\in N_{1}$ and
$d_{\mathcal{F} - X_{1-i}, W}^i(v)\geq \eta n^{k}$ for $v\in (N
\setminus e_0)\cap [n]$. Thus there exists an edge $e_{u+1}$ in $\mathcal{F} - \cup_{i=0}^u e_i$ such that $\{v_{u+1}\}= e_{u+1}\cap N$, $|e_{u+1}\cap W|\equiv 0\pmod 2$ when $u+1\leq r$, $|e_{u+1}\cap W|\equiv 1\pmod 2$ when $r < u+1  \le s$, and $|e_{u+1}\cap W|\equiv i\pmod 2$ and $e_j\cap X\subseteq X_i$ when $s+1 < u+1  \le t$.
Continuing this process for at most $|N\setminus e_0|$ steps, we obtain the
desired matching $M$. \qed

Let ${\cal H}(W,U;r)$ denote the balanced $(1,k)$-partite $(k+1)$-graph
with partition classes $X, [n]$ and  a partition $W,U$  of
$[n]$ such that for every $x\in X$, $N_{{\cal H}(W,U;r)}(x)=
\{e\in {[n]\choose k}: |e\cap W|=r\}$.
Now we show that if all the vertices of a $(1,k)$-partite
$(k+1)$-graph ${\cal F}$ are good with respect to ${\cal H}(W,U;r)$,
then there exists a perfect matching in ${\cal F}$
consisting of edges intersecting $W$ exactly $r$ times.

\begin{lemma}\label{good-PM}
Let $k,n, r$ be integers such that $k \ge 3$, $0\leq r\leq k$,  $n
\equiv 0 \pmod k$, and $1/n\ll 1/k$.
Let $\mathcal{F}$ be a balanced $(1,k)$-partite $(k+1)$-graph with partition classes
$X, [n] $, and let $W, U$ be a partition of $[n]$ with $|W| = rn/k$.
Suppose all vertices in $\mathcal{F}$ are
 {$\gamma$-good} with respect to  ${\cal H}(W,U;r)$.
Then there exists a perfect matching $M$ in $\mathcal{F}$ such that
$|e\cap W|=r$ for all $e\in M$.
\end{lemma}

\pf By symmetry between $W$ and $U$, we may assume $r\ge k/2$.
Let $M$ be a maximum matching in $\mathcal{F}$ such that, for every
$e\in M$, $|e\cap W|=r$.
Let $W_0 := W \setminus V(M)$ and $U_0 := U \setminus V(M)$.
Since $r \ge k/2$, $|W_0|\ge |U_0|$.

Suppose $|M| < \frac{n}{4k}$. Then  $|W_0| \ge \frac{n}{4}$. By
maximality of $M$, for $v\in W_0$, we have
$$|N_{H(W,U;r)}(v) \setminus N_{\mathcal{F}}(v)|
\ge |X \backslash V(M)| {|W_0| \choose r-1} {|U_0| \choose k-r}. $$

Thus, if $r=k$ then
$$|N_{H(W,U;r)}(v) \setminus N_{\mathcal{F}}(v)|\ge \frac{3n}{4k}  {|W_0| \choose k-1}\ge \gamma  (n+n/k)^{k}.$$

Now  suppose $r \le k-1$. Then
$$|U_0| =\frac{k-r}{r}|W_0|\ge
    \frac{k-r}{r} (n/4)\ge \frac{n}{4(k-1)}.$$
So we have
\begin{align*}
|N_{H(W,U;r)}(v) \setminus  N_{\mathcal{F}}(v)|
&\ge \frac{3n}{4k} (|W_0| - r + 2)^{r-1} (|U_0| - k +r + 1)^{k-r}/k! \\
&\ge \gamma (n+n/k)^{k}  
\end{align*}
contradicting the fact that $v$ is  {$\gamma$-good} with $\mathcal{H}(W,U;r)$.

Now, suppose for a contradiction that $M$ is not a perfect matching.
There exist $x_{k+1} \in X\backslash V(M)$, distinct $v_{k+1,1},...,v_{k+1,r} \in W_0$,
and distinct $v_{k+1,r+1}, ... ,v_{k+1,k} \in U_0$.
Let $\{e_1,e_2,...,e_{k}\}\in {M\choose k}$
and write $e_i := \{ x_i, v_{i,1}, ..., v_{i,k} \}$, such that,  for
$i \in [k]$, $x_i \in X$, $v_{i,j} \in W$ for $j \in [r]$, and
$v_{i,j} \in U$ for $j \in [k]\setminus [r]$.
For $j \in [k+1]$, let $f_j := \{x_j, v_{j+1,1}, v_{j+2,2}, ..., v_{j+1+k,k}\}$, with the addition in the subscripts modulo $k+1$ (except we write $k+1$ for $0$).
Note that $f_1,...,f_{k+1}$ are pairwise disjoint and $|f_j\cap W|=r$
for $j\in [k+1]$.

If $f_j \in E(\mathcal{F})$ for all $j \in [k+1]$, then $M' := (M \cup
\{f_1,...,f_{k+1}\})\setminus \{e_1,...,e_{k}\}$ is matching in $\mathcal{F}$, contradicting the maximality of $|M|$.
Hence, $f_j \not\in E(\mathcal{F})$ for some $j \in [k+1]$.
Note that there are ${|M| \choose k}$ choices of $\{e_1,...,e_k\} \subseteq M$.
Thus we have
\begin{align*}
&\Big|\{e \in E(H(W,U;r)) \setminus E({\mathcal{F}}) : |e \cap \{v_{k+1,i}: i \in [k]\}| = 1\}\Big| \\
&\ge {|M| \choose k}  \\
&\ge \left(\frac{n}{4k} - k + 1\right)^k / k! \\
&\ge \frac{1}{k!(5k)^k} n^k \text{\quad { (since $n \ge 20k^2$) } }\\
&> \gamma(k+1) (n+n/k)^k 
\end{align*}
This implies that there exists $x \in \{x_{k+1},v_{k+1,1},\ldots,v_{k+1,k}\}$ such that
$$| N_{H(W,U;r)}(x) \setminus N_{\mathcal{F}}(x) | > \gamma (n+n/k)^{k}.$$
That is,  $x$ is not {$\gamma$-good} with respect to $H(W,U;r)$, a contradiction.
\qed

\medskip

After obtaining the matching $M$ in $\mathcal{F}$ from Lemma \ref{good-PM}, we need to find a perfect matching in $\mathcal{F} - V(M)$ to
conclude the proof of Theorem~\ref{general} for the case when all
vertice of  $\mathcal{F} - V(M)$ are good. The following lemma
serves this purpose.

\begin{lemma} \label{close-final}
Let $k,n$ be integers such that $k \ge 3$, $n \equiv 0 \pmod k$, and
$n\gg 1/\varepsilon$.   Let $\mathcal{F}$ be a balanced $(1,k)$-partite $(k+1)$-graph with partition
classes $X, [n]$. Let $W',U'$ be a
partition of $[n]$ with $\min\{|U'|,|W'|\}\geq 1.1n/k+k$, and let
$i\in \{0,1\}$ such that
\begin{itemize}
\item[(i)] if $i=0$ then $|W' |\equiv 0\pmod 2$,
\item[(ii)] if $i=1$ then $|W'| \equiv n/k \pmod 2$,
\item[(iii)] if $k-i\equiv 0 \pmod 2$ then $|U'|\equiv 0\pmod 2$, and
\item[(iv)] if $k-i\equiv 1 \pmod 2$ then $|U'| \equiv n/k \pmod 2$.
\end{itemize}
If all vertices of $\mathcal{F}$ are  {$\gamma$-good} with
respect to $\mathcal{H}_{n,k}^i(W',U')$, then there is a perfect
matching in $\mathcal{F}$.
\end{lemma}

\pf First, suppose there exists an edge $e_1\in E({\cal F})$ such that $|W'\setminus
e_1|=r_1x+r_2y$ and $|U'\setminus e_1|=(k-r_1)x+(k-r_2)y$, where $x,y,
r_1,r_2$ are integers satisfying $x,y > 20k^2$ and $0\le r_1,r_2\le
k$.  We partition $W'\setminus e_1$ to $W_1,W_2$ such that $|W_1|=r_1x$ and
$|W_2|=r_2y$, partition $U'\setminus e_1$ to $U_1,U_2$ such
that $|U_1|=(k-r_1)x$ and $|U_2|=(k-r_2)y$, and partition $X\setminus
e_1$ to $X_1,X_2$ such that $|X_1|=x$ and $|X_2|=y$.
By assumption, for $i\in \{0,1\}$, $0\le r_i\le k$.
Hence, by Lemma~\ref{good-PM}, ${\cal F}[X_i\cup W_i\cup U_i]$ has a perfect
matching, say  $M_i$, consisting of edges containing exactly $r_i$
vertices from $W_i$. Now $M_1\cup M_2$ (when $e_1=\emptyset$) or
$M_1\cup M_2\cup \{e_1\}$ gives the desired matching.

Therefore, it suffices to prove the existence of such $e_1$. This is
done in the following four cases.

\medskip

{\it Case} 1.   $i=0$ and $k-i=k$ is even.

Then (i) and (iii) hold. So $|W'|$ and $|U'|$ are both
even. Then, since all vertices of $\mathcal{F}$ are  $\gamma$-good with respect to $\mathcal{H}_{n,k}^0(W',U')$,
there is an edge $e_1\in \mathcal{F}$ such that $|e_1 \cap W'| \equiv
|W'| \pmod k$. Thus, $|W'\setminus e_1|\equiv 0 \pmod k$.

Since $n\equiv 0 \pmod k$, we have $|e_1\cap
U'|=k-|e_1\cap W'|\equiv
k-|W'|\pmod k$; hence, $|e_1\cap U'|\equiv |U'|\pmod k$. So
$|U'\setminus e_1|\equiv 0 \pmod k$.

Thus, we may take $r_1=k$ and $r_2=0$.
Note that $x = |W'\setminus e_1|/k > n/k^2 > 20k^2$ and $y = |U'\setminus e_1|/k > n/k^2 > 20k^2$.

\medskip

{\it Case} 2. $i=0$ and $k-i=k$ is odd.

Then (i) and (iv) hold. So $|W'|$ is even and $|U'| \equiv n/k \pmod
2$. Since all vertices of  $\mathcal{F}$ are $\gamma$-good with respect to
$\mathcal{H}_{n,k}^0(W',U')$,
there is an edge $e_1\in \mathcal{F}$ such that $|e_1 \cap W'| \equiv
|W'| \pmod {k-1}$.

Write $|W'\setminus e_1|=(k-1)x$ for some
integer $x$.
Then $| U'\setminus e_1|-x=n-x-|W'|-|e_1\cap U'|=n-x-|W'\backslash e_1|-|W'\cap e_1|+|e_1\cap U'|=n-k(x+1)\equiv 0
\pmod k$.
 So we take $r_1=k-1$ and $r_2=0$.

Note that $x = |W'\setminus e_1|/(k-1) > n/k^2 > 20k^2$
 and 
\begin{align*}
 y &= (| U'\setminus e_1|-x)/k \\
 &= (| U'\setminus e_1|-|W'\setminus e_1|/(k-1))/k\\
 &=  (| U'\setminus e_1|-(n-k-| U'\setminus e_1|)/(k-1))/k\\
 &= (| U'\setminus e_1| - n/k + 1) / (k-1)\\
 &\ge 0.1n/k^2 > 20k^2.
\end{align*}

\medskip

{\it Case} 3.  $i=1$ and $k-i=k-1$ is odd.

Then (ii) and (iv) hold. So $|W'| \equiv n/k \pmod 2$ and $|U'| \equiv n/k \pmod 2$.
Without loss generality, suppose that $|W'|\geq |U'|$.
Since all vertices  of $\mathcal{F}$ are {$\gamma$-good} with respect to $\mathcal{H}_{n,k}^1(W',U')$,
there is an edge $e_1\in E(\mathcal{F})$ such that
$|W' \setminus e_1| \equiv n/k-1\pmod {k-2}$.
So  $|e_1 \cap W'|\equiv |W'| - n/k + 1 \pmod {k-2}$.

Let $|W'\setminus e_1|=(k-1)x+y$ and $x+y=n/k-1$; then $x+(k-1)y=n-k-|W'\setminus e_1|=|U'\setminus
e_1|$. Moreover, $x=(|W'\setminus e_1|-n/k+1)/(k-2)$ and $y=(|U'\setminus e_1|-n/k+1)/(k-2)$.  So we may set $r_1=k-1$ and $r_2=1$.

Note that $x=(|W'\setminus e_1|-n/k+1)/(k-2) > 0.1n/k^2 > 20k^2$ and
$y=(|W'\setminus e_1|-n/k+1)/(k-2) > 0.1n/k^2 > 20k^2$.

\medskip

{\it Case} 4.   $i=1$ and $k-i=k-1$ is even.

Then (ii) and (iii) hold.
So $|W'| \equiv n/k \pmod 2$ and $|U'|$ is even. Since all  vertices
in $\mathcal{F}$ are  $\gamma$-good with respect to
$\mathcal{H}_{n,k}^1(W',U')$,
there is an edge $e_1\in E(\mathcal{F})$ such that $|e_1 \cap U'|
\equiv |U'| \pmod {k-1}$.

Write $|U'\setminus e_1|=(k-1)y$ for some
integer $y$. 
Then 
\begin{align*}
|W'\setminus e_1|-y&= n-y -|U'|-|e_1\cap W'|\\
&=n-y -|U'\backslash e_1|-(|e_1\cap U'|+|e_1\cap W'|)\\
&=n-y-(k-1)y-k\\
&=n-k(y+1)\equiv
0\pmod k.
\end{align*} So let $r_1=k$ and $r_2=1$.

Note that
\begin{align*}
x & = (| W'\setminus e_1|-y)/k \\
&= (| W'\setminus e_1|-|U'\setminus e_1|/(k-1))/k\\
  & =  (| W'\setminus e_1|-(n-k-| W'\setminus e_1|)/(k-1))/k\\
  & = (| W'\setminus e_1| - n/k + 1) / (k-1)\\
  &\ge 0.1n/k^2 \\
  & > 20k^2.
\end{align*}
Moreover $y = |U'\setminus e_1|/(k-1) > n/k^2 > 20k^2$. \qed

\medskip

We are ready to prove the main result in this section.

\begin{lemma} \label{close}
Let $k,n$ be integers with $k \ge 3$, $n \equiv 0 \pmod k$, and $n\gg
1/\varepsilon$.
Let $\mathcal{F}$ be a balanced $(1,k)$-partite $(k+1)$-graph with partition
classes  $X, [n]$, such that $\delta_{k-1}(N_{\mathcal{F}}(x)) >
t(n,k)$ for all $x \in X$.  Suppose $[n]$ has a partition
$W, U$ with $|W| = n/2 \pm o(n)$ and $|U|=n/2 \pm o(n)$ such that $\mathcal{F}$ is
strongly $\varepsilon$-close to ${\cal H}_{n,k}^0(W,U;m)$ for some $m\in
[n/k]$.
Then  $\mathcal{F}$ admits a perfect matching.
\end{lemma}

{
\pf Since $\mathcal{F}$ is
strongly $\varepsilon$-close to ${\cal H}_{n,k}^0(W,U;m)$ for some $m\in
[n/k]$, the number of $\varepsilon^{2/3}$-bad vertices in $\mathcal{F}$ is at most $(k+1)\varepsilon^{1/3}(n/k+n)$ $\varepsilon^{2/3}$ (see Lemma \ref{bad-vertices}). Let $B$ denote the set of $\varepsilon^{2/3}$-bad vertices in $\mathcal{F}$. Write $X_2=X\cap B$ and $N=[n]\cap B$.
For $i\in \{0,1\}$,
\[
X_i=\{x\in X\backslash X_2\ |\ N_{{\cal H}_{n,k}^0(W,U;m)}(x)\cong H_{n,k}^i(W,U)\}.
\]
Write $m_i:=|X_i|$ for $i\in \{0,1,2\}$. For $i\in \{0,1\}$ and $x\in X_i$,  $x$ is $\varepsilon^{2/3}$-good; hence,
\[
|N_{{\cal H}^0(W,U;m)}(x)\setminus N_{\mathcal{F}}(x)|\leq \varepsilon^{2/3}(n+n/k)^k\leq (1+1/k)^k\varepsilon^{2/3} n^k.
\]
So $N_{\mathcal{F}}(x)$ is strongly $(1+1/k)^k\varepsilon^{2/3}$-close to $H_{n,k}^i(W,U)$.
 
Recall the definition of constants $\eta = 1/(4k!)$,
$c = 1/(8(k+1)!)$,
$\varepsilon = 1/(80^kk^{k-5}(k!)^k((k+1)!)^k)^{3/2}$ and $\gamma = \frac{\varepsilon^{2/3} 2^k}{c^k k^k}$. 
Let $\mathcal{F}_i:=(\mathcal{F}-X_2)-X_{1-i}$ for $i\in \{0,1\}$.
Define a partition $W_0,U_0$ of $[n]$ as follows: If $|X_0|\geq |X_1|$ then let
$$W_0=( W \setminus  \{v \in W \cap N: d_{\mathcal{F}_0, W}^0(v) \leq \eta n^{k} \} ) \cup \{v \in U \cap N: d_{\mathcal{F}_0, W}^1(v) >\eta
n^{k}  \}$$
and
$$U_0=( U \setminus \{v \in U \cap N: d_{\mathcal{F}_0, W}^1(v) > \eta
n^{k} \} ) \cup \{v \in W \cap N: d_{\mathcal{F}_0, W}^0(v) \leq \eta
n^{k} \}. $$
If $|X_0|<|X_1|$ then let
$$W_0=( W \setminus  \{v \in W \cap N: d_{\mathcal{F}_1, W}^1(v) \leq
\eta n^{k} \} ) \cup \{v \in U \cap N: d_{\mathcal{F}_1, W}^0(v) >\eta
n^{k}  \}$$
and
$$U_0=( U \setminus \{v \in U \cap N: d_{\mathcal{F}_1, W}^0(v) > \eta n^{k}
\} ) \cup \{v \in W \cap N: d_{\mathcal{F}_1, W}^1(v) \leq \eta n^{k}
\}.$$

 Let $X_{21} := \{x\in X_{2}: d_{\mathcal{F}, W_0}^0(x) \geq \eta n^{k}
\} $ and $X_{22}=X_2\backslash X_{21}$.
We apply Lemma \ref{parity} to $\mathcal{F}$ with $X' = X_{1} \cup
X_{22}$ and $i=0$ (if $|X_0| \ge |X_1|$), or $X' = X_{0} \cup
X_{21}$ and $i=1$ (when $|X_0|<|X_1|$). So there exists $e_0$ such that
$e_0=\emptyset$ or $e_0\in E(\mathcal{F})$ satisfying (i)(ii)(iii)(iv) of Lemma \ref{parity}.

Define
\begin{equation*}
N' :=
\begin{cases}
B \cup X_1, & \text{if $|X_1| \le  cn$} \\
B \cup X_0, & \text{if $|X_0|\le cn$} \\
B, & \text{otherwise}.
\end{cases}
\end{equation*} Then $|N'|\leq 2cn$.
Next, we apply Lemma \ref{match-bad-vertices} to $\mathcal{F}$ with
$W_0$ and $U_0$ as partition of $[n]$ and with $N'$ as the $N$ in Lemma~\ref{match-bad-vertices}.
Note that if $|X_0| \ge |X_1|$, we set $i=0$ and have $d_{\mathcal{F} - X_{1}, W_0}^0(v)\geq
\eta n^{k}$ for all $v\in [n]$, and  that if $|X_0|< |X_1|$, we set $i=1$ and have
$d_{\mathcal{F} - X_{0}, W_0}^1(v)\geq \eta n^{k}$ for all $v\in [n]$.
Hence, there is a matching $M_1$ in
$\mathcal{F}$ with $e_0\in M_1$ (when $e_0\ne \emptyset$), and 
$|V(M_1)| \le (k+1)2cn$
such that
\begin{itemize}
\item[(i)] $N' \subseteq V(M_1)$;
\item[(ii)] $|W_0 \setminus V(M_1)| \equiv  |X_{1} \backslash V(M_1)| \pmod 2$.
\end{itemize}

Let $W_1:= W_0\setminus V(M_1)$ and $U_1:=U_0 \setminus V(M_1)$.
Let $X_0':=X_0\setminus V(M_1)$ and $X_1'=X_1\setminus  V(M_1)$.
Note that  $|X_i'| > cn/2$ or $X_i' = \emptyset$ for $i \in \{0,1\}$.
If $|X_0|<cn$, then $X_0'=\emptyset$ and let $W_{11}:=W_1$, $U_{11}:=U_1$ and $W_{10}=U_{10}=\emptyset$; if $|X_1|<cn$, then $X_1'=\emptyset$ and let $W_{10}:=W_1$, $U_{10}:=U_1$ and $W_{11}=U_{11}=\emptyset$. 
Otherwise, let $W_{10} \cup W_{11}$ be a partition of $W_1$ and $U_{10} \cup U_{11}$ be a partition of $U_1$ such that for $i\in \{0,1\}$,
\begin{itemize}
\item [(iii)] $|W_{10}|\equiv 0\pmod 2$, $|W_{11}|\equiv |X_{1}'|\pmod 2$, $|W_{1i}\cup U_{1i}|=k|X_{i}'|$, and if $X_{i}' \ne \emptyset$ then $\min\{|W_{1i}|, |U_{1i}|\} \ge 1.1|X_{i}'| + k$. 
\end{itemize}
Note that the existence of these partitions is guaranteed by (ii).

Let $\mathcal{F}_i':=\mathcal{F}[X_i'\cup W_{1i}\cup U_{1i}]$ for $i\in \{0,1\}$. 
Since every vertex in $V(\mathcal{F}) \backslash B$ is $\varepsilon^{2/3}$-good with respect to $\mathcal{H}_{n,k}^0(W,U;m)$,
every vertex of $\mathcal{F}_i'$ is $\gamma$-good  with
respect to
$\mathcal{H}_{k|X_i'|,k}^i(W_{1i},U_{1i})$. 
For, otherwise, without loss of generality, suppose $v \in V(\mathcal{F}_i')$ is $\gamma$-bad with respect to  $\mathcal{H}_{k|X_{i}'|,k}^i(W_{1i},U_{1i})$.
Then
\begin{align*}
|N_{\mathcal{H}_{n,k}(W,U;m)}(v) \setminus N_{\mathcal{F}}(v)|
& \ge |N_{\mathcal{H}_{k|X_i'|,k}^i(W_{1i},U_{1i})}(v) \setminus
  N_{\mathcal{F}_{i}'}(v)| \\
&> \gamma ((k+1)|X_i'|)^k  \\
&\geq \frac{\varepsilon^{2/3} 2^k}{c^k k^k} \left(\frac{(k+1)cn}{2}\right)^k  \\
&= \varepsilon^{2/3} (n+n/k)^k,
\end{align*}
contradicting that $v$ in $\mathcal{F}$ is $\varepsilon^{2/3}$-good on $\mathcal{H}_{n,k}(W,U;m)$.

By Lemma \ref{close-final},
$\mathcal{F}_i'$ contains a perfect matching $M_{2i}$ for $i=0,1$ (and let $M_{2i} = \emptyset$ if $\mathcal{F}_i'$ is empty).
So $M_1\cup M_{20}\cup M_{21}$ is a perfect matching in $\mathcal{F}$.
\qed

\medskip

}



\section{Absorbing devices for perfect matchings}

We need the following lemma from \cite{RRS09}. For subsets $N_1,
\ldots, N_k$ of the vertex sets of a $k$-graph $H$, let $E_H(N_1,
\ldots, N_k):= \{(v_1, \ldots, v_k): v_i\in N_i \mbox{ for } i\in [k]
\mbox{ and } \{v_1,\ldots, v_k\}\in E(H)\}$. Let $e_H(N_1, \ldots, N_k):=|E_H(N_1, \ldots, N_k)|$.
Given a $k$-graph $H$, let $\overline{H}$ denote the $k$-graph with vertex set $V(H)$ and edge set
\[
E(\overline{H})=\left\{e\in {V(H)\choose k} : e\notin E(H)\right\}.
\]
The following lemma is Claim 5.1 in \cite{RRS09}.



\begin{lemma}[R\"odl, Ruci\'nski, and Szemer\'edi \cite{RRS09}] \label{absorb-2}
Let $n,k$ be two integers such that $n\gg k\geq 3$ and $n\equiv 0\pmod k$. Let $H$ be a $k$-graph on $n$ vertices. If $\delta_{k-1}(H)\geq (1/2-1/\log n)n$ and $H$ is not weakly $\varepsilon$-close to $H(n,k)$ or $\overline{H(n,k)}$, then at least one of the following holds.
\begin{enumerate}
  \item [(i)] For all $N_1,\ldots, N_k\subseteq V(H)$ with $|N_i|\geq
    (1/2-1/\log n)n$, we have $e_H(N_1,\ldots, N_k)\geq n^k/\log^3 n$.

  \item [(ii)] $|\{(v_1,\ldots,v_{k-1})\in {V(H)\choose {k-1}}:
    d_H(\{v_1,\ldots,v_{k-1}\})>(1/2+2/\log n)n\}|\ge  n^{k-1}/\log n.$
\end{enumerate}
\end{lemma}

Next we define two types of absorbing devices for a given balanced set $S$ of $k+1$ vertices.
Both are  $(k+1)$-matchings. The vertices of each devices together with $S$ induce a $(k+1)$-graph with a perfect matching.

\smallskip

{\bf Absorbing device I}:
Let $\mathcal{H}$ be a $(1,k)$-partite $(k+1)$-graph with partition classes $X,V$.
Given a balanced set $S =\{x_0,v_1,\ldots, v_k\}$ with $x_0 \in X$ and
$v_i \in V$ for $i \in [k]$, a $(k+1)$-matching
$\{e_1,\ldots,e_{k},g\}$ in $\mathcal{H}$ is said to be {\it
  $S$-absorbing} if $\mathcal{H}$ has a $(k + 2)$-matching $\{e_1',\ldots,e_{k}',f,g\}$ such that
\begin{itemize}
  \item [(i)] $e_i'\cap e_j=\emptyset $ for all $i\neq j$,
  \item [(ii)] $e_i'\setminus e_i=\{v_i\}$ and $|e_i\setminus e_i'|=1$
    for $i \in [k]$, and
  \item [(iii)] $f=\{x_0\}\cup \bigcup_{i\in [k]}(e_i\setminus e_i')$.
\end{itemize}
Note that the inclusion of the edge $g$ is only for later convenience.

\smallskip

{\bf Absorbing device II}:
Let $\mathcal{H}$ be a $(1,k)$-partite $(k+1)$-graph with partition classes $X,V$.
Given a balanced set $S =\{x_0,v_1,\ldots, v_k\}$ with $x_0 \in X$ and
$v_i \in V$ for $i \in [k]$, a $(k+1)$-matching
$\{e_1,\ldots,e_{k+1}\}$ in $\mathcal{H}$ is said to be {\it
  $S$-absorbing} if $\mathcal{H}$ has a $(k + 2)$-matching $\{e_1',\ldots,e_{k+1}',f\}$ such that
\begin{itemize}
\item [(i)] $e_i\cap e_j'=\emptyset$ for $i,j\in [k]$, where $i\neq j$,

\item [(ii)] $e_i' \setminus e_i= \{v_i\}$   and $|e_i\setminus e_i'|=1$ for $i \in [k]$,

\item [(iii)] $e_{k+1}'\cap e_k=e_k\setminus e_k', |e_{k+1}\setminus e_{k+1}'|=1$;

\item [(iv)] $f=\{x_0\}\cup \bigcup_{i\in [k+1]\setminus \{k\}}(e_i\setminus e_i')$
\end{itemize}

Next, we show that if (i) or (ii) of Lemma \ref{absorb-2} holds for $\mathcal{F}(n,k)$,
then for each balanced set $S$ there are many $S$-absorbing devices in
$\mathcal{F}(n,k)$.

\begin{lemma}\label{absorb-count}
Let $\mathcal{F}$ be a $(1,k)$-partite $(k+1)$-graph with partition
classes $X,[n]$ such that $\delta_{k-1}(N_{\mathcal{F}}(x)) \geq (1/2-1/\log
n)n$ for each $x\in X$. Let \\$R:=\{x\in X : N_{\mathcal{F}}(x) \mbox{  is not weakly
  $\varepsilon$-close to } H(n,k) \mbox{ or } \overline{H(n,k)}\}$.
Let $S =\{x_0,v_1,\ldots, v_k\}$  be a balanced $(k+1)$-set such that $x_0\in R$.
\begin{enumerate}
\item[(i)] If (i) of Lemma \ref{absorb-2} holds for $N_{\mathcal{F}}(x_0)$, then the number of $S$-absorbing devices I in $\mathcal{F}$ is $\Omega(n^{(k+1)^2}/\log^3 n)$.
\item[(ii)] If $|R|=o(n)$ and $|\{x\in R\setminus \{x_0\}: \mbox{  (ii) of Lemma
    \ref{absorb-2} holds for } N_{\mathcal{F}}(x)\}|\ge n/\log n$, then the number of $S$-absorbing devices II in $\mathcal{F}$ is $\Omega(n^{(k+1)^2}/\log^4 n)$.
\end{enumerate}
\end{lemma}

\pf
First, suppose (i) of Lemma~\ref{absorb-2} holds for
$N_{\mathcal{F}}(x_0)$. Since $\delta_{k-1}(N_{\mathcal{F}}(x_0)) \ge
(1/2-1/\log n)n$,  there are $\Omega(n^{k})$ sets $B_i$, for each for $i\in [k]$, such that $B_i\cup \{v_i\}\in \mathcal{F}$. Consequently, there are
$\Omega(n^{k^2})$ choices of (pairwise disjoint) such sets
$B_1,\ldots, B_k$. Let $e_i'= B_i\cup \{v_i\}$ for $i\in [k]$.

Since $|N_{\mathcal{F}}(B_i)|\geq (1/2-1/\log n)n$ and  (i) of Lemma~\ref{absorb-2} holds for
$N_{\mathcal{F}}(x_0)$, we have
\[
e_{\mathcal{F}}(x_0,N_{\mathcal{F}}(B_1),\ldots, N_{\mathcal{F}}(B_k))\geq n^{k}/\log^3 n.
\]
So, there are at least $n^{k}/\log^3 n$ choices of edges
$f=\{x_0,u_1,\ldots,u_k\}$ such that $e_i:=B_i\cup \{u_i\}\in
E(\mathcal{F})$ for $i\in [k]$. Moreover, there are $\Omega(n^{k+1})$ choices $g$ such that $g\in E((\mathcal{F}-S)-\cup_{i=1}^k e_i)$. Hence, altogether there are
\[
(n^{k}/\log^3 n)\Omega(n^{k^2})\Omega(n^{k+1})=\Omega(n^{(k+1)^2}/\log^3 n)
\]
choices of $S$-absorbing $(k+1)$-matchings $\{e_1',
\ldots,e_k',f,g\}$.

Next we show (ii). As in the argument for (i), since $|R|=o(n)$, there are
$\Omega(n^{k^2})$ choices of (pairwise disjoint) sets $B_1,\ldots,B_k$
such that $R \cap (\cup_{i=1}^k B_i) = \emptyset$ and $e_i'=B_i\cup
\{v_i\}\in E(\mathcal{F})$ for $i\in [k]$.

For $i \in [k-1]$, we choose $u_i\in  N_{\mathcal{F}} (B_i)$, each in
at least $(1/2-1/\log n)n$ ways, and let $e_i=B_i\cup \{u_i\}$.
  Let  $ y \in R \setminus \{x_0\}$ such that (ii) of Lemma
  \ref{absorb-2} holds for $N_{\mathcal{F}}(y)$. By assumption, there
  are at least $n/\log n $ such $y$.
  We select a $(k-1)$-element sequence of vertices, say $T$, such that
  $d_{N_{\mathcal{F}(y)}}(T)\ge (1/2+2/\log n)n$ and $T$ is disjoint from $S\cup B_k\cup \bigcup_{i=1}^{k-1} e_{i}$.
Let $B_{k+1}=T\cup \{y\}$.  Then we may pick $u_k$, $u_{k+1}$  such
that $B_k\cup \{u_k\}, B_{k+1}\cup \{u_k\}, B_{k+1}\cup \{u_{k+1}\}$,
$\{u_0,u_1,\ldots,u_{k-1},u_{k+1}\}\in E(\mathcal{F})$. Note that
                                $d_{\mathcal{F}}(B_{k+1})\geq (1/2+2/\log n)n$.
We have $|N_{\mathcal{F}}(B_{k+1})\cap N_{\mathcal{F}}(B_{k})|\geq n/\log n$ and $|N_{\mathcal{F}}(B_{k+1})\cap N_{\mathcal{F}}(\{x_0,u_1\ldots,u_{k-1}\})\ |\geq n/\log n$. Thus there are at least $n/\log n$ choices for each of $x_k,x_{k+1}$.
  By Lemma \ref{absorb-2} (ii), there are at least $n^{k-1}/\log n$ choices for $T$. Summarizing, we have chosen $B_1,\ldots, B_k, B_{k+1}, u_1,\ldots,u_{k-1}$, $T$ forming an $S$-absorbing $(k+1)$-matching, in
  \[
  \Omega(n^{k^2})(n^{k-1}/\log n) n^{k-1} (n/\log n)^3)=\Omega(n^{(k+1)^2}/\log^4 n)
  \]
ways. \qed

\smallskip

{\bf Absorbing device III}:
Let $\mathcal{H}$ be a $(1,k)$-partite $(k+1)$-graph with partition classes $X,V$.
Given a balanced set $S =\{x_0,v_1,\ldots, v_k\}$ with $x_0 \in X$ and
$v_i \in V$ for $i \in [k]$, a $(k+1)$-matching
$\{e_1,\ldots,e_{k+1}\}$ in $\mathcal{H}$ is said to be {\it
  $S$-absorbing} if $\mathcal{H}$ has a $(k + 2)$-matching $\{e_1',\ldots,e_{k+1}',f\}$ such that
\begin{itemize}
\item [(i)] $e_i\cap e_j'=\emptyset$ for $i,j\in [k]$, where $i\neq j$,

\item [(ii)] $e_i' \setminus e_i= \{v_i\}$   and $|e_i\setminus e_i'|=1$ for $i \in [k]$,

\item [(iii)] $e_{k+1}'\setminus e_{k+1}=\{x_0\}, |e_{k+1}\setminus e_{k+1}'|=1$,

\item [(iv)] $f=\bigcup_{i\in [k+1]}(e_i\setminus e_i')$.
\end{itemize}


In the following lemma, we show that any $(1,k)$-partite
$(k+1)$-graph is close to a $(k+1)$-graph with same partition classes, or contains an absorbing matching.

\begin{lemma}\label{close-absorbing}
Let $0<\varepsilon\ll \beta\ll 1$. Let $\mathcal{H}$ be a
balanced $(1,k)$-partite $(k+1)$-graph with partition classes $X, [n]$.
Let $R:=\{x\in X: N_{\mathcal{H}}(x) \mbox{ is not weakly
  $\varepsilon$-close to } H_{n,k}^0 \mbox{ or }H_{n,k}^1\}$.  Let
$x_0\in X\setminus R$ such that
$N_{\mathcal{H}}(x_0)$ is strongly $\varepsilon$-close to $
H_{n,k}^{i}(W,U)$ for some $i\in \{0,1\}$ and for some partition
$(W,U)$ of $[n]$ with $|W|=(1/2\pm o(1))n=|U|$.
Suppose $|R|\leq n/\log n$ and $\delta_{k-1}(N_{\mathcal{H}}(x))\geq
(1/2-1/\log n)n$ for every $x\in X$.
If $\mathcal{H}$ is not strongly $\beta$-close to $\mathcal{H}^0_{n,k}(W,U; m)$
for any integer $0 \le m \le n/k$, then for any balanced $(k+1)$-set $S$ with $x_0\in S$,
there are at least $\Omega(n^{(k+1)^2}/\log ^4n)$ $S$-absorbing
devices I or III.
\end{lemma}

\pf First, let $\alpha$ be a number with $\varepsilon\ll \alpha\ll \beta$.  We prove

\textbf{Claim 1.~}  Let $x\in X$ and let $W_x, U_x$ be a partition of $[n]$ such that $|W_x|=(1/2\pm o(1))n$ and $|U_x|=(1/2\pm o(1))n$ and $N_{\mathcal{H}}(x)$
is strongly $\varepsilon$-close to $H_{n,k}^j(W_x,U_x)$ for some $j\in
\{0,1\}$.  For any  $N_1,\ldots, N_k\subseteq [n]$ with $|N_i|\geq(1/2-1/\log n)n$ for $i \in [k]$,
 if $e_{N_{\mathcal{H}}(x)}(N_1,N_2,\ldots,N_k)\leq n^{k}/\log n$ then,  for $j\in [k]$,  $|U_x\setminus N_j|\leq \alpha n$ or $|W_x\setminus N_j|\leq \alpha n$.

Otherwise, suppose, without loss of generality,   $|U_x\setminus N_1|\geq \alpha n$ and $|W_x\setminus N_1|\geq \alpha n$.  Let $D$ denote the set
of $\sqrt{\varepsilon}$-bad vertices on $H_{n.,k}^j(W_x,U_x)$; then $|D|\le
k\sqrt{\varepsilon}n$.
Then, since $|N_1|\geq (1/2-1/\log n)n$, $| N_1\setminus (U_x\cup D)|\geq \alpha n/2$ and $|N_1\setminus (W_x\cup D)|\geq \alpha n/2$.
Note for any $z \in N_1\setminus (W_x\cup D)$ and $y\in N_1\setminus
(U_x\cup D)$, since $\varepsilon \ll \alpha$ and both $z$ and $y$ are $\sqrt{\varepsilon}$-good with respect to $H_{n,k}^j(W_x,U_x)$, we have
$|N_{\mathcal{H}}(\{x,z\})\cup N_{\mathcal{H}}(\{x,y\})|\geq {n-1\choose k-1}-\alpha n^{k-1}$.
Hence, either $e_{N_{\mathcal{H}}(x)}(\{z\},N_2,\ldots,N_k)\geq
\frac{n^{k-1}}{2^kk!}-\alpha n^{k-1}\geq \frac{n^{k-1}}{2^{k+1}k!}$ or
$e_{N_{\mathcal{H}}(x)}(y,N_2,\ldots,N_k)\geq
\frac{n^{k-1}}{2^kk!}-\alpha n^{k-1}\geq
\frac{n^{k-1}}{2^{k+1}k!}$. So by symmetry, we may assume that there exists $N_1'\subseteq N_1\setminus (W_x\cup  D)$ such that 
 $|N_1'|\geq |N_1\setminus (W_x\cup  D)|/2$ and for all $z\in N_1'$
\[
e_{N_{\mathcal{H}}(x)}(\{z\},N_2,\ldots,N_k)\geq \frac{n^{k-1}}{2^{k+1}k!},
\]
Hence,
\[
e_{N_{\mathcal{H}}(x)}(N_1,N_2,\ldots,N_k)\geq \frac{\alpha n^{k}}{2^{k+3}k!}>n^k/\log n,
\]
contradicting the assumption of Claim 1 and completing its proof.

\medskip

 Let $A=\{x\in X\setminus R : |N_{\mathcal{H}}(x)\cap N_{\mathcal{H}}(x_0)|\geq
\alpha n^k\}$ and let $B=X\setminus (A\cup R)$.
Let $x\in B$. Then
\begin{align*}
 |N_{\mathcal{H}}(x)\cap E(H_{n,k}^i(W,U))|&\leq |N_{\mathcal{H}}(x)\cap N_{\mathcal{H}}(x_0)|+|E(H_{n,k}^i(W,U))\setminus N_{\mathcal{H}}(x_0)|
 \leq  \alpha n^k+\varepsilon n^k.
\end{align*}
Hence, for any $x\in B$, we have
\begin{align*}
|E(\overline{H_{n,k}^i(W,U)})\setminus E(N_{\mathcal{H}}(x))|
&= e(\overline{H_{n,k}^i(W,U))})-(|N_{\mathcal{H}}(x)|-|N_{\mathcal{H}}(x)\cap E(H_{n,k}^i(W,U))|)\\
&\leq e(\overline{H_{n,k}^i})-(e(H_{n,k}^i)-\alpha n^k-\varepsilon n^{k}) +  |N_{\mathcal{H}}(x)\cap E(H_{n,k}^i(W,U))|\\
&\leq \varepsilon n^k+\alpha n^k+\varepsilon n^k+\alpha n^k+\varepsilon n^k\leq 3\alpha n^k.
\end{align*}
Thus, we have

\medskip

{\bf Claim} 2. For any $x\in B$, $N_{\mathcal{H}}(x)$ is strongly $(3\alpha)$-close to $\overline{H_{n,k}^i(W,U)}$,
which is strongly   $3\alpha$-close to $H_{n,k}^{1-i}(W,U)$. So 

\medskip

\textbf{Claim} 3. $|B|\leq (1-\beta)n/k$; so $|A|\ge \beta n/k -n/\log n$.

For, otherwise, $|B| > (1-\beta)n/k$. Then $|A|+|R|<\beta n/k$ and
\begin{align*}
|E(\mathcal{H}_{n,k}^{1-i}(W,U; |B|)) \setminus E(\mathcal{H})|&=\sum_{x\in B}|N_{\mathcal{H}_{n,k}^{1-i}(W,U;|B|)}(x)\setminus N_{\mathcal{H}}(x)|+(|A|+|R|){n\choose k}\\
&<  \frac{n}{k}4\alpha n^k+\beta n/k {n\choose k} \leq \beta (n+n/k)^{k+1},
\end{align*}
a contradiction since $\mathcal{H}$ is not strongly $\beta$-close to $\mathcal{H}_{n,k}^{1-i}(W,U; |B|)$.

\medskip


\textbf{Claim} 4. Let ${\cal H}'={\cal H}[A\cup [n]]$. Then, for any  $N_1,\ldots, N_k\subseteq [n]$ with $|N_i|\geq (1/2-1/\log n)n$, either $e_{\mathcal{H}'}(\{x_0\},N_1,\ldots,N_k)\geq n^{k}/\log n$ or $e_{\mathcal{H}'}(A,N_1,\ldots,N_k)\geq n^{k+1}/\log^3n $.

Suppose on the contrary that there exist $N_1,\ldots, N_k\subseteq [n]$ with
$|N_i|\ge  (1/2-1/\log n)n$, such that
\begin{align}\label{eq:1}
e_{\mathcal{H}'}(\{x_0\},N_1,\ldots,N_k)< n^{k}/\log n
\end{align}
and
\begin{align}\label{ab-eq:2}
e_{\mathcal{H}'}(A,N_1,\ldots,N_k)< n^{k+1}/\log^3n.
\end{align}
Since $N_{\mathcal{H}}(x_0)$ is strongly $\varepsilon$-close to
$H_{n,k}^i(W,U)$, the number of
$\sqrt{\varepsilon}$-bad vertices in $N_{{\cal H}}(x)$ with respect
to  $H_{n,k}^i(W,U)$ is at most $k\sqrt{\varepsilon}n$.

By Claim 1 and (1), $|U\setminus N_i|\leq \alpha n$ or $|W\setminus N_i|\leq \alpha n$ for each $i \in [k]$.
Let $A_1:=\{x\in A: e_{N_{\mathcal{H}}(x)}(N_1,\ldots, N_k)\geq
n^{k}/\log n\}$ and $A_2=A\setminus A_1$. By (\ref{ab-eq:2}), $|A_1|\leq n/\log^2 n$.

Let $y\in A_2$. So $y \not\in R$ and, thus, $N_{\mathcal{H}}(y)$ is weakly $\varepsilon$-close to
$H_{n,k}^0$ or $H_{n,k}^1$.
Let  $W_y, U_y$ be a partition of $[n]$ such that $N_{\mathcal{H}}(y)$
is strongly $\varepsilon$-close to $H_{n,k}^j(W_y,U_y)$ for some $j\in
\{0,1\}$.
By Claim 1, we have for all $i\in[k]$, $|U_y\setminus N_i|\leq \alpha n$ or $|W_y\setminus N_i|\leq \alpha n$. Hence
\begin{align}\label{ab-eq:3}
|U_y\setminus U|\leq |U_y\setminus N_i|+|N_i\setminus U|\leq 2\alpha n
\end{align}
or
\begin{align}\label{ab-eq:4}
|U_y\setminus W|\leq |U_y\setminus N_i|+|N_i\setminus W|\leq 2\alpha n
\end{align}

 We claim that $N_{\mathcal{H}}(y)$ is strongly $(5\alpha)$-close to $H_{n,k}^j(W,U)$ if inequality (\ref{ab-eq:3}) holds.
\begin{align*}
|E(H_{n,k}^j(W,U))\setminus N_{\mathcal{H}}(y) |&\leq | E(H_{n,k}^j(W_y,U_y))\setminus N_{\mathcal{H}}(y)|+|  E(H_{n,k}^j(W,U)) \setminus E(H_{n,k}^j(W_y,U_y))|\\
&\leq \varepsilon n^k+4\alpha n^k \leq 5\alpha n^k,
\end{align*}
If inequality (\ref{ab-eq:4}) holds then
\begin{align*}
| E(H_{n,k}^j(U,W))\setminus N_{\mathcal{H}}(y)|&\leq | E(H_{n,k}^j(W_y,U_y))\setminus  N_{\mathcal{H}}(y)|+| E(H_{n,k}^j(U,W))\setminus E(H_{n,k}^j(W_y,U_y))| \\
&\leq \varepsilon n^k+4\alpha n^k \leq 5\alpha n^k.
\end{align*}
Thus  $N_{\mathcal{H}}(y)$ is strongly $(5\alpha)$-close to $H_{n,k}^j(U,W)$. 

Hence by Claim 2, for all $x\in X\setminus (A_1\cup R)$, $N_{\mathcal{H}}(x)$ is strongly $(5\alpha)$-close to $H_{n,k}^j(W,U)$ for some $j\in \{0,1\}$.
 For $j\in \{0,1\}$, let $X_j=\{x \in X: \mbox{$N_{\mathcal{H}}(x)$
   is strongly $(5\alpha)$-close to $H_{n,k}^{j}(W,U)$}\}$.
Since $|A_1| \le n/\log^2 n$,
and $|R|\le n/\log n$, we have $|X\setminus (X_0\cup X_1)|=|A_1\cup R|
\leq2 n/\log n$.
Hence
\begin{align*}
|E(\mathcal{H}_{n,k}^{0}(W,U;|X_0|))\setminus E(\mathcal{H})|&\leq \sum_{j=0}^1\sum_{x\in X_j}|N_{\mathcal{H}_{n,k}^{0}(W,U;|X_0|)}(x)\setminus N_{\mathcal{H}}(x)|+(|A_1|+|R|)n^k\\
&\leq \frac{n}{k}5\alpha n^k+2n^{k+1}/\log n\leq \beta (n+n/k)^{k+1}.
\end{align*}
So  $\mathcal{H}$ is strongly $\beta$-close
$\mathcal{H}_{n,k}^{0}(W,U;|X_0|)$, a contradiction.
This concludes the proof of Claim 4.

\medskip

Let $S=\{x_0,v_1,\ldots,v_k\}$ with $v_1, \ldots, v_k\in [n]$ all
distinct. For each $i\in [k]$, there are $\Omega(n^{k})$ sets $B_i$
such that $e_i' =B_i\cup \{v_i\}\in E(\mathcal{H})$, and there are
$\Omega(n^{k^2})$ choices of (pairwise disjoint) such sets $B_1,\ldots,
B_k$. 
For $i\in [k]$, we have $N_{\mathcal{H}'}(B_i)\geq (1/2-1/\log n)n$ by assumption. So   there are  $\Omega(n^{k^2})$ choices of (disjoint) sets $B_1,\ldots, B_k$ such that by Claim 4, 
one of the following two inequalities holds:
\begin{align}
&e_{\mathcal{H}'}(\{x_0\},N_{\mathcal{H}'}(B_1),\ldots,N_{\mathcal{H}'}(B_k))\geq n^{k}/\log n\label{cl-ab-eq1}\\
&e_{\mathcal{H}'}(A,N_{\mathcal{H}'}(B_1),\ldots,N_{\mathcal{H}'}(B_k))\geq n^{k+1}/\log^3n \label{cl-ab-eq2}.
\end{align}

First, assume (\ref{cl-ab-eq1}) holds. Then there at least $n^k/\log
n$ choices $u_1,\ldots, u_k$ such that $e_i=B_i\cup \{u_i\}\in
E(\mathcal{H'})$ for $i\in [k]$ and $\{x_0,u_1,\ldots,u_k\}\in
E(\mathcal{H'})$. Moreover, there are $\Omega(n^{k+1})$ choices of
$g\in E(\mathcal{H'})$ such that $g$ is disjoint from we $\cup_{i=1}^k
e_i$. Thus the number of  $S$-absorbing devices I is at least
\[
\Omega(n^{k^2}) \Omega(n^k/\log
n)\Omega(n^{k+1})=\Omega(n^{(k+1)^2}/\log n).
\]

Now assume (\ref{cl-ab-eq2}) holds. Then there are $n^{k+1}/\log^3n$
choices $\{y,u_1,\ldots,u_k\}$ such that $y\in A$ and,
for $i\in [k]$, $B_i\cup \{u_i\}\in E(\mathcal{H'})$.
By definition of $A$, there are
$\Omega(n^{k})$ choices $B_{k+1}$ such that $B_{k+1}\in
N_{\mathcal{H}}(y)\cap N_{\mathcal{H}}(x_0)$. Let $e_{k+1}=B_{k+1}\cup
\{y\}$. So there are at least
$\Omega(n^{k+1}/\log^3n)\Omega(n^{k})\Omega(n^{k^2})=\Omega(n^{(k+1)^2}/\log
^3n)$ different choices of $\{e_1,\ldots,e_{k+1}\}$ such that
$\{e_1,\ldots,e_{k+1}\}$ is an $S$-absorbing device  III. \qed

\medskip

To  prove another absorbing lemma,
we need to use Chernoff bounds, see \cite{AS08}.

\begin{lemma}
\label{chernoff}
Suppose $X_1, ..., X_n$ are independent random variables taking values in $\{0, 1\}$. Let $X$ denote their sum and $\mu = \mathbb{E}[X]$ denote the expected value of $X$. Then for any $0 < \delta \leq 1$,
$$\mathbb{P}[X \leq (1-\delta) \mu] < e^{-\frac{\delta^2 \mu}{2}}.$$
\end{lemma}

\begin{lemma}\label{Absorb-lem}
Let $k \geq 3$ be an integer and let $0<\varepsilon\ll \beta\ll 1$.
There exists $n_2 > 0$ such that the following holds for any integer $n > n_2$.
Let $\mathcal{H}$ be a balanced
$(1,k)$-partite $(k+1)$-graph with partition
classes $X,[n]$ such that  and $\delta_{k-1}(N_{\mathcal{H}}(x))>
t(n,k)$ for all $x\in X$.
 Let $R:=\{x\in X: N_{\mathcal{H}}(x) \mbox{  is not weakly
   $\varepsilon$-close to } H_{n,k}^0 \mbox{  or }H_{n,k}^1\}$, and
 let $x_0\in X$.
   Suppose one of the following three conditions holds for every $S\in {V(H)\choose k+1}$ with $S\cap X=\{x_0\}$:
\begin{itemize}
\item[$(i)$] (i) of Lemma \ref{absorb-2} holds for
  $N_{\mathcal{H}}(x_0)$, and $\mathcal{H}$ has
  $\Omega(n^{(k+1)^2}/\log^4 n)$ $S$-absorbing  devices I.
\item[$(ii)$]  $|\{x\in R \setminus  \{x_0\}:
  \mbox{   (ii) of Lemma \ref{absorb-2} holds for }
  N_{\mathcal{H}}(x)\}|\ge n/\log n$, and $\mathcal{H}$ has
  $\Omega(n^{(k+1)^2}/\log^4 n)$ $S$-absorbing devices II.
\item [$(iii)$] $|R| \le n/\log n$, $x_0 \notin R$, and $N_{\mathcal{H}}(x_0)$ is
  strongly $\varepsilon$-close to $H_{n,k}^0(W,U)$ or
  $H_{n,k}^1(W,U)$ for some partition $(W,U)$ of $[n]$ with $|W| = (1/2 + o(1))n = |U|$,
$\mathcal{H}$ is not strongly $\beta$-close to ${\cal
  H}_{n,k}^0(W,U;m)$ for any $m\in [n/k]$, and   $\mathcal{H}$ has
$\Omega(n^{(k+1)^2}/\log^4 n)$ $S$-absorbing devices I or III.
\end{itemize}
Then there exists a matching $M'$ in $\mathcal{H}$ such that
$|M'|=O(\log^6 n)$ and, for each $S\in {V(H)\choose k+1}$ with $S\cap
X=\{x_0\}$, $M'$ contains an $S$-absorbing $(k+1)$-matching.
\end{lemma}

\pf For each balanced  $(k+1)$-set $S \subseteq V(\mathcal{H})$ with $x_0\in S$,
let $\Gamma(S)$ be the collection of $S$-absorbing $(k+1)$-matchings.
 Then by Lemmas \ref{absorb-count} and \ref{close-absorbing},  $|\Gamma(S)|  =\Omega(n^{(k+1)^2}/\log^4 n)$.
 So we may choose  constant $\alpha:=\alpha(k) > 0$ such that
 \[|\Gamma(S)| \geq \alpha (k^2+k)!{n/k\choose k+1}{n\choose k(k+1)}/((k!)^{k+1}\log^4 n).\]

Let $\mathcal{M}$ be the  family
obtained by choosing  a sequence of balanced $(k+1)$-sets $(S_1,\ldots, S_{k+1})$ independently with probability
\[
p=\frac{(k!)^{k+1}\log^6 n }{{n/k\choose k+1}{n\choose k(k+1)}(k^2+k)!}.
\]
Note that $p<1$ as we can choose $n_2$ large enough. Then
\[
\mathbb{E}(|\mathcal{M}|)=p  {n/k\choose k+1}{n\choose k(k+1)}(k^2+k)!/(k!)^{k+1}=O(\log^6 n),
\]
and, for  $(k+1)$-set $S\subseteq V(\mathcal{H})$ with $\{x_0\}=S\cap X$,
\[
\mathbb{E}(|\mathcal{M} \cap \Gamma(S)|) \geq p   \alpha (k^2+k)!{n/k\choose k+1}{n\choose k(k+1)}/((k!)^{k+1}\log^4 n) = \alpha \log^2 n.
\]

By Lemma~\ref{chernoff} and by choosing $n_2$ large enough, we have,
for $n>n_2$ and for each $S\in {V(H)\choose k+1}$ with $S\cap X=\{x_0\}$,
$$
\mathbb{P}[|\mathcal{M}| > 2\alpha \log^6 n] = \mathbb{P}[|\mathcal{M}| > 2 \mathbb{E}(|\mathcal{M}|)] \leq e^{-\mathbb{E}(|\mathcal{M}|)/3} = e^{- (\log^6 n)/3}.
$$
So with probability at least $1-o(1)$
\begin{align}\label{prop-F}
|\mathcal{M}|\leq 2 \alpha\log^6 n.
\end{align}
Again by Lemma~\ref{chernoff} and  by choosing $n_2$ large enough, we
have, for $n>n_2$ and for ,
\begin{align*}
\mathbb{P}[|\mathcal{M} \cap \Gamma(S)| \leq (\alpha \log^2 n)/2 ] &\leq  \mathbb{P}[|\mathcal{M} \cap \Gamma(S)| \leq  \mathbb{E}(|\mathcal{M} \cap \Gamma(S)|)/2 ] \\
& \leq  e^{-\mathbb{E}(|\mathcal{M} \cap \Gamma(S)|)/8}\\
&\leq  e^{- (\alpha\log^2 n)/8}.
\end{align*}
So by union bound and by choosing $n_2$  large, we have for $n>n_2$,
\begin{align*}
&\mathbb{P}[\exists S \in {V(H)\choose k+1} \mbox{ with } S\cap X=\{x_0\} : |\mathcal{M} \cap \Gamma(S)|
 \leq (\alpha \log^2 n)/2 ] \\
& \leq {n\choose k} e^{- (\alpha\log^2 n)/8} \\
& = 2n^{k - (\alpha\log n)/8} < 1/10.
\end{align*}
Thus,  with probability at least $9/10$, for all $S \in {V(H)\choose
  k+1}$ with $S\cap X=\{x_0\}$, we have
\begin{align}\label{prop-F-S}
|\mathcal{M} \cap \Gamma(S)| \geq (\alpha\log^2 n)/2 > 1.
\end{align}

Furthermore, the expected number of pairs of sequences $(S_1, \ldots, S_{k+1}), (T_1, \ldots, T_{k+1})\in {\cal M}$
satisfying $(\cup_{i \in [k]} S_i) \cap (\cup_{i \in [k]} T_i) \neq \emptyset$ is at most
\begin{align*}
&\left(\frac{(k^2+k)!}{(k!)^k}\right)^2{n/k\choose k+1}{n\choose k^2+k} \\
&\left({k+1\choose 1}{n/k\choose k}{n\choose k^2+k}+ {k^2+k\choose 1}{n/k\choose k}{n\choose k^2+k-1}\right)p^2<1/2.
\end{align*}
Thus, with probability at least $1/2$ (by Markov's inequality), for all distinct $(S_1,\ldots, S_{k+1})\in {\cal M}$ and
$(T_1,\ldots, T_{k+1})\in \mathcal{M}$,
\begin{align}\label{F-cap}
\mbox{$\bigcup_{i\in[k+1]} S_i$ and $\bigcup_{i\in [k+1]} T_i$ are disjoint.}
\end{align}

Hence, with positive probability, $\mathcal{M}$ satisfies (\ref{prop-F}),  (\ref{prop-F-S}), and (\ref{F-cap}). So we may assume that ${\cal M}$
satisfies (\ref{prop-F}),  (\ref{prop-F-S}), and (\ref{F-cap}).
Let $M$ be the union of ${\cal M}\cap \Gamma(S)$ for all $S\in
{V(H)\choose k+1}$ with $S\cap X=\{x_0\}$.
Then $M$ is the desired matching. \qed

\section{Absorbing devices for near perfect matchings}

Let $H$ be a $(1,k)$-partite $(k+1)$-graph with  partition classes $Q,
[n]$. 
For a set $S\in {V(H)\choose k+2}$ with $|S\cap Q|=1$, an edge $e\in E(H)$ is said to
be {\it $S$-absorbing} if $H[e\cup S]$ has a matching  of size $2$.


\begin{lemma}
\label{absorbing-counting}
Let $k\ge 3$, $0 < c < 1$, and $n$ be a sufficiently large integer, and
let $H$ be a $(1,k)$-partite $(k+1)$-graph with partition classes
$Q,[n]$ such that $|Q|\leq (n-1)/k$. If
$\delta_{k-1}(N_H(v))\geq cn$ for all $v\in Q$, then, for any $S\in {V(H)\choose k+2}$
with $|S\cap Q|=1$, $H$  has at least $c^4n^{k+1}/2$ $S$-absorbing edges.
\end{lemma}

\pf Let $S=\{v\}\cup B$, where $v\in Q$ and
$B=\{b_1,\ldots,b_k,b_{k+1}\}\in {[n]\choose k+1}$.  Let
$B'=B\setminus \{b_k,b_{k+1}\}$.
Since $\delta_{k-1}(N_H(v))\geq cn$, we have at least $cn-2$ choices
for $q\in [n]\setminus B$ such
that $B'\cup \{v,q\} \in E(H)$. For a $(k-2)$-set $A$ of
$V(H)\setminus (S\cup \{q\})$ with $|A\cap Q|=1$, since $d_{H}(A\cup
\{b_k,b_{k+1}\})\geq cn$, there are at least $cn-2k$ choices $q'\in
[n]\setminus (S\cup A\cup \{q\})$ such that $\{q'\}\cup \{b_{k},b_{k+1}\}\cup A\in E(H)$.
Since $d_{H}(A\cup \{q,q'\})\geq cn$, there are at least $cn-2k$
choices $q''$ such that $A\cup \{q,q',q''\}\in E(H)$. Clearly, any such $\{q,q',q''\}\cup A$ is $S$-absorbing.

Note that there are $(|Q|-1){n-k-2\choose k-3}$ different choices for set $A$.  Hence the number of $S$-absorbing edges is at least
\[
(|Q|-1){n-k-2\choose k-3}(cn-2)(cn-2k)(cn-2k)\geq c^4n^{k+1}/2.
\]
This completes the proof. \qed

\medskip

Analogous to the absorbing results in \cite{RRS09,Han15, Lu16},
we prove the existence of small absorbing matching
for $(1,k)$-partite $(k+1)$-graphs with $|Q|\ge n/k-k^2$.
\begin{lemma}
\label{absorbing-near}
For any constant $c>  0$, there exists an integer  $n_0 > 0$ with the following
holds for every integer $n\ge n_0$: Let $H$ be a $(1,k)$-partite $(k+1)$-graph with partition classes
$Q, [n]$ such that  $n/k-k^2\leq |Q|\leq (n-1)/k$
and $\delta(N_H(v)) \geq cn$ for all $v\in S$.
Then there is a matching $M$ in $H$ such that $|M|\le
(32(k+3)/c^4)\log n$ and, for each $S\in {V(H)\choose k+2}$
with $|S\cap Q|=1$, $M$ contains at least  $4(k+3)\log n$ $S$-absorbing edges.
\end{lemma}

\pf  Let $C = 32(k+3)/c^4$.
Let $M$ be the family obtained by choosing each edge independently with probability
$p = (C/2) n^{-(k+1)} \log n $. Thus $\mathbb{E}[|M'|] =|E(H)|p\le n^{k+1} p = (C/2) \log n$.

The number of  intersecting pairs of edges in $E(H)$ is at most
$|Q|{n\choose k}^2+|Q|^2n{n-1\choose k-1}^2\leq n^{2k+1}$; so
the expected number of  intersecting pairs of edges in $M'$ is at most $$n^{2k+1} p^2 \leq  C^2 \log^2 n/(4n) = o(1).$$
By Markov's inequality, with probability strictly larger than $1/3$, $M$ is a matching of size at most $C \log n$.

For a set $S\in {V(H)\choose k+2}$ with $|S\cap Q|=1$, let $X_S$ denote the number of $S$-absorbing edges in $M$.
Then by Lemma \ref{absorbing-counting}, we have
$$\mathbb{E}[X_S] \geq p c^4n^{k+1}/2  = 8(k+3)\log n.$$
By Lemma~\ref{chernoff},
$$\mathbb{P}[X_S \leq \ \mathbb{E}[X_S] /2] \leq \exp(-
\mathbb{E}[X_S]/8) = \exp(-(k+2)\log n) = n^{-(k+2)}.$$

Note that there are at most $|Q|{n\choose k}< n^{k}/k$ sets $S\in
{V(H)\choose k+1}$ with $|S\cap Q|=1$.  It follows from union bound that,  with probability strictly larger than $1/4$,
 $X_S\ge \mathbb{E}[X_S]/2 \geq 4(k+3)\log n$  for all $S\in {V(H)\choose k+1}$ with $|S\cap Q|=1$. Thus, the desired $M$
 exists.
\qed

\section{Hypergraphs not close to extremal configurations} 

In this section, we prove Theorem~\ref{main} for  hypergraphs that are not close to extremal configurations.
For this, we need a result on almost perfect matchings in $(1,k)$-partite $(k+1)$-graphs.

\begin{lemma}\label{large-mat}
Let  $k,n$ be positive integers with $k\geq 3$. Let $H$ be a  $(1,k)$-partite $(k+1)$-graph with vertex
partition classes $Q, [n]$, where $k+1 \leq |Q|\leq n/k$. If $\delta_{k-1}(N_H(v))>n/k$ for every $v\in Q$, then $H$ has a matching covering all but at most $k-1$ vertices of $Q$.

\end{lemma}

\pf  Let $M$ be a maximum matching in $H$. We may
assume $|Q\setminus V(M)| \ge k$;
for, otherwise, $M$ gives the desired matching. Note  $$|[n]\setminus
V(M)|=n-k|M|\ge k|Q|-k|M|\ge k|Q\setminus V(M)|\geq k^2.$$
So there exist $k$ pair-disjoint $k$-sets in $[n]\backslash V(M)$, say $S_1,\ldots, S_k$ such that $|S_i\cap Q|=1$ for $i\in [k]$.
By maximality  of $M$,  $N_H(S_i)\subseteq V(M)$ for $i\in [k]$.

We claim that $\sum_{i=1}^k |N_H(S_i)\cap e|\le k$ for all $e\in M$. For
otherwise, there exist $e\in M$ and distinct $u,v\in e$ such that
$u\in N_H(S_i)\cap e$ and $v\in N_H(S_j)\cap e$. Now $(M\setminus
\{e\})\cup \{S_i\cup \{u\}, S_j\cup \{v\}\}$ is a matching in $H$,
contradicting the maximality of $M$.

Since  $N_H(S_i)\subseteq V(M)$,
$$
\sum_{i=1}^k |N_H(S_i)|
=\sum_{e\in M} \sum_{i=1}^k |N_H(S_i)\cap e|
\leq k|M|
< n.$$
However, since $\delta_{k-1}(N_H(v)) >n/k$ for all $v\in Q$, we have $\sum_{i=1}^k
|N_H(S_i)| \ge  k \delta_{k-1}(N_H(v)) > n$,
 a contradiction.  \qed

\begin{lemma}\label{near-pm}
Let $k,n$ be integers with $k\geq 3$. Let $H$ be a $(1,k)$-partite
$(k+1)$-graph with partition classes $Q, [n]$, where $k+1 \leq |Q|\leq
(n-1)/k$. If $\delta_{k-1}(N_H(v))\geq n/k+32k^5(k+3)\log n$ for every $v\in Q$, then $H$ has a matching covering $Q$.
\end{lemma}

\pf Let $c=1/k$. Since
 $\delta_{k-1}(N_H(v))\geq n/k+32k^5(k+3)\log n$ for every $v\in Q$,
 by Lemma \ref{absorbing-near}, $H$ has a matching $M$ of size at most
 $32k^4(k+3)\log n$ such that for any set $S\in {V(H)\choose k+2}$
 with $|S\cap Q|=1$, the number of $S$-absorbing edges in $M$ is at
 least $k+1$. Let $H'=H-V(M)$, with partition classes $Q\setminus
 V(M), [n]\setminus V(M)$.
For every $v\in Q\setminus V(M)$, we have
 \[
 \delta_{k-1}(N_{H'}(v))\geq \delta_{k-1}(N_H(v))-k|M|>n/k.
 \]
Thus by Lemma \ref{large-mat}, $H'$ has a matching $M'$ covering all
but at most $k$ vertices of $Q\setminus V(M)$.

Let $M_0:=M$ and $M_0'=M'$. If $M_0\cup M_0'$ covers $Q$, then we are
done. Otherwise, there exists $S_1\in {V(H)\setminus V(M_0\cup
M_0')\choose k+2}$ with $|S_1\cap Q|=1$.
Recall $M_0$ contains an  $S_1$-absorbing edge, say $e_0$; so
$H[S_1\cup e_0]$ contains a matching $X_0$ of size 2. Then
$M_1':=M_0'\cup X_0$ is a matching in $H$. Let $M_1:=M_0\setminus \{e_0\}$.
If $M_1\cup M_1'$ covers $Q$, then we are done. Otherwise, since
$|Q\setminus V(M_1\cup M_1')|\le k$ and $M$ has at least $k+1$
$S$-absorbing matching for each $S\in {V(H)\setminus V(M_1\cup
M_1')\choose k+2}$ with $|S\cap Q|=1$, we may repeat
the above procedure.    Thus, we obtain a  maximal sequence of pairs
of matchings $M_0, M_0', M_1, M_1', \ldots, M_{t}, M_t'$,
a sequence of $(k+2)$-sets $S_1,\ldots, S_{t}$ with $|S_i\cap Q|=1$
for $i\in [t]$, and a sequence of matchings $X_0, X_1, \ldots, X_{t}$. Now $M_t\cup M_t'$ is a matching of $H$ covering $Q$.  \qed

\begin{lemma} \label{not-close}
Let $\varepsilon > 0$  be a constant and $k,n$ be integers with  $k \ge 3$ and $n \equiv 0 \pmod k$
such that $0<1/n\ll \varepsilon\ll 1/k$.
Let $\mathcal{F}$ be a balanced $(1,k)$-partite $(k+1)$-graph with partition
classes $X$ and $[n]$.  Let $R:=\{x\in X:
N_{\mathcal{F}}(x) \mbox{  is not weakly $\varepsilon$-close to }
H_{n,k}^0 \mbox{ or } H_{n,k}^1\}$.
Suppose
\begin{itemize}
\item $\delta_{k-1}(N_{\mathcal{F}}(x)) >  t(n,k)$ for all $x \in X$;
\item $|R| > n/\log n$,  or $|R|\leq n/\log n$ and $\mathcal{F}$ is
  not strongly $\varepsilon$-close to ${\cal H}_{n,k}^0(W,U;m)$ for
  any $m\in [n/k]$ and for any partition $W,U$ of $[n]$ with $|W| = n/2 \pm o(n)$ and $|U|=n/2 \pm o(n)$.
\end{itemize}
Then  $\mathcal{F}$ admits a perfect matching.
\end{lemma}

\pf Note that the conclusion of Lemma~\ref{Absorb-lem} holds. We define $x_0 \in X$ as follows:
\begin{itemize}
\item[(i)] If $|R| > n/\log n$ then choose $x_0\in R$ such that,  whenever
  possible, (i)
  of Lemma \ref{absorb-2} holds for  $N_{\mathcal{F}}(x_0)$.

\item [(ii)] If $|R| \le n/\log n$ then
let $x_0 \in X\setminus R$ and let $W, U$ be a partition of $[n]$ with $|W| = n/2 \pm o(n)=|U|$ such that $N_{\mathcal{F}}(x_0)$ is strongly $\varepsilon$-close to $H_{n,k}^0(W,U)$ or $H_{n,k}^1(W,U)$.
\end{itemize}

By Lemma \ref{Absorb-lem}, there exists a matching $M$ in ${\cal F}$
with $|M|=O(\log^6n )$ such that for any balanced $(k+1)$-set
$S\subseteq V({\cal F})$ containing $x_0$, $\mathcal{F}[S\cup V(M)]$ has a perfect matching.

 Let $\mathcal{F}'=\mathcal{F}- (V(M) \cup \{x_0\})$. Then  $k+1\le |V(\mathcal{F}')\cap X|=|V(\mathcal{F}')\cap [n]|/k - 1$. Moreover, for every $v\in X\cap V(\mathcal{F}')$, we have
 \[
 \delta_{k-1}(N_{\mathcal{F}'}(v))\geq \delta_{k-1}(N_{\mathcal{F}}(v)) - |V(M)|\geq n/2-k-k|M|>n/2-k-k\log^6 n>n/k+\log^2 n.
 \]
 Thus by Lemma \ref{near-pm}, $\mathcal{F}'$ has a matching $M'$ covering $ X\backslash( V(M) \cup \{x_0\})$.
Now it is easy to see that $S:=V(\mathcal{F})-V(M\cup M')$ is a balanced $(k+1)$-set with $\{x_0\}= S\cap X$.
Note that $\mathcal{F}[S\cup V(M)]$ has a perfect matching $M''$.
Therefore, $M'\cup M''$ is a perfect matching of $\mathcal{F}$.
\qed

\section{Concluding remarks}

First, we point out that Theorem~\ref{general} follows immediately
from Lemmas \ref{close} and \ref{not-close}, which in turn implies
Theorem~\ref{main}.

Thus, we proved a rainbow version of the result of R\"odl,
Ruci\'nski, and Szemer\'edi \cite{RRS09} that determines the
co-degree threshold function for the existence of a perfect matching
in a $k$-graph.

There are many results on various Dirac type conditions for the
existence of a matching of certain size. One can ask questions about
whether similar rainbow versions hold for those results. Our method of
converting the rainbow matching problem to a matching problem for a
special class of hypergraphs provides a way for establishing the rainbow
versions by using existing tools for matching problems.

We list some results that their rainbow version may be studied using
our approach.  R\"odl,
Ruci\'nski, and Szemer\'edi \cite{RRS09} proved that, for $n \not\equiv 0\pmod k$, the minimum co-degree threshold that ensures a matching $M$ in a $k$-graph $H$
with $|V(M)| \geq |V(H)| - k$ is between $\lfloor n/k\rfloor$ and $n/k+O(\log n)$, and  conjectured that this threshold function is $\lfloor n/k\rfloor$. This conjecture  was proved recently
by Han \cite{Han15}.
Treglown and Zhao \cite{TZ12, TZ13} determined the minimum $l$-degree threshold for perfect matchings in $k$-graphs for $k/2 \leq l \leq k-1$. Bollob\'{a}s, Daykin,  and Erd\H{o}s \cite{BDE} considered the minimum vertex degree for the appearance of matchings of certain size.
They proved that for integer $k\geq 2$, if $H$ is a $k$-graph of order $n\geq 2k^3(m-1)$ and
$\delta_1(H)>{{n-1}\choose {k-1}}-{{n-m}\choose {k-1}},$
then $H$ has a matching of size at least $m$.
The bound on $n$ is improved to $n \ge 3k^2 m$ by Huang and Zhao \cite{HZ17} recently.
For $3$-graphs, K\"{u}hn, Osthus, and Treglown \cite{KOT} and, independently, Khan \cite{Khan13} determined the minimum vertex degree threshold for perfect matchings in $3$-graphs, which improves an earlier result by H\`an, Person, and Schacht \cite{HPS}.
The minimum vertex degree threshold for perfect matchings in $4$-graphs is obtained by Khan in \cite{Khan16}.

\end{document}